\documentclass[10pt, conference, letterpaper]{IEEEtran}
\makeatletter
\def\ps@headings{%
\def\@oddhead{\mbox{}\scriptsize\rightmark \hfil \thepage}%
\def\@evenhead{\scriptsize\thepage \hfil \leftmark\mbox{}}%
\def\@oddfoot{}%
\def\@evenfoot{}}
\makeatother
\usepackage[pdftex]{graphicx}
\usepackage{amsmath}
\usepackage{amssymb}
\usepackage{dblfloatfix}
\usepackage{array}
\usepackage[svgnames]{xcolor} 
\usepackage{mdwmath}
\usepackage{subfigure}
 \usepackage{bm}
\usepackage{tipa}
 \usepackage{hyperref}
 \usepackage{framed}
 \usepackage{booktabs}
 \usepackage{caption}
 \usepackage{url}
\usepackage{algorithm}
\usepackage{algorithmic}
\usepackage{epsfig}
\usepackage{overpic}
\usepackage{epstopdf}
\usepackage{amsfonts,setspace,amsthm,braket}

\graphicspath{{../pdf/}{../jpeg/}}
\DeclareGraphicsExtensions{.pdf,.jpeg,.png}
\DeclareMathOperator*{\argmin}{arg\,min}

\usepackage{dsfont}
\newtheorem{theorem}{Theorem}
\newtheorem{lemma}{Lemma}

\newtheorem{definition}{Definition}

\allowdisplaybreaks

\newcommand\EatDot[1]{}
\newcommand{\edit}[1]{{ #1}}

\floatname{algorithm}{Algorithm}
\IEEEoverridecommandlockouts

\begin{document}
%\title{Throughput-Optimal Broadcast Algorithm in Wireless Networks} %original
%\title{Dynamic Broadcast in Wireless Networks} %georgios
%\title{Optimal Distributed Broadcast Algorithm for Wireless Networks} %abhishek
%\title{Dynamic Algorithm for Broadcast in a DAG} %original
\title{Network Utility Maximization with Heterogeneous Traffic Flows} %chihping
%\title{A Dynamic Broadcast Algorithm in Directed Acyclic Graphs} %chihping

% Commenting the author section below out

%
\author{
%\phantom{
\IEEEauthorblockN{Abhishek Sinha, Eytan Modiano}
\IEEEauthorblockA{Laboratory for Information and Decision Systems, Massachusetts Institute of Technology, Cambridge, MA 02139}
%\IEEEauthorblockA{\IEEEauthorrefmark{2}Mathematical and Algorithmic Sciences Lab France Research Center, Huawei Technologies Co., Ltd.\\}
%\IEEEauthorblockA{\IEEEauthorrefmark{3}Qualcomm Research, San Diego, CA\\
Email: sinhaa@mit.edu,
modiano@mit.edu

%\thanks{ This work was supported by the NSF Grants CNS-1217048 and CNS-1524317.}
%\thanks{\IEEEauthorreqfmark{2} The work of G. Paschos was done while he was at MIT and affiliated with CERTH-ITI, and it was supported in part by the WiNC project of the Action: Supporting Postdoctoral Researchers, funded by national and Community funds (European Social Fund).}
%\thanks{\IEEEauthorrefmark{3} The work of C.p.Li was done when he was a Postdoctoral scholar at LIDS, MIT.}

}
\maketitle
\begin{abstract} \label{abstract_section}
%  What we do: we extend the virtual queue framework of  \cite{UMW} to NUM problems, and allow dynamic variation of the network. Second, we connect the VQ framework, which was obtained by relaxing the precedence constraints, to the dual formulation of the NUM problem. We show that, the VQ appears naturally when we consider a path-based formulation of the NUM problem. The virtual queues correspond to the dual variables, when the dual problem is solved by subgradient algorithm with a constant step-size of one. We also derive explicit form of the corresponding dual objective function for logarithmic and $\alpha$-fair utilities. 
%We consider the Network Utility Maximization (NUM) problem for wireless networks in the presence of arbitrary types of flows - including unicast, broadcast, multicast and anycast traffic. Building upon the recent framework of a universal control policy (UMW), we design a utility optimal cross-layer admission control, routing and scheduling policy, called \textbf{UMW+}. The \textbf{UMW+} policy takes packet level actions, based on a \emph{precedence-relaxed} virtual network. We establish a precise one-to-one correspondence between the dynamics of the  virtual queues under the UMW+ policy, and the dynamics of the dual variables of an associated offline NUM program, under a subgradient algorithm. Extensive numerical simulation results  highlight the efficacy of the proposed policy. 
We consider the Network Utility Maximization (NUM) problem for wireless networks in the presence of arbitrary types of flows, including unicast, broadcast, multicast, and anycast traffic. Building upon the recent framework of a universal control policy (UMW), we design a utility optimal cross-layer admission control, routing and scheduling policy, called UMW+. The UMW+ policy takes packet level actions based on a \emph{precedence-relaxed} virtual network.  Using Lyapunov optimization techniques, we show that UMW+ maximizes network utility, while simultaneously keeping the physical queues in the network stable.  Extensive simulation results validate the performance of UMW+; demonstrating both optimal utility performance and bounded average queue occupancy. Moreover, we establish a precise one-to-one correspondence between the dynamics of the virtual queues under the UMW+ policy, and the dynamics of the dual variables of an associated offline NUM program, under a subgradient algorithm.    This correspondence sheds further insight into our understanding of UMW+. 
\end{abstract}

\section{Introduction} \label{intro_section}
%%\textbf{TBD}.
%Outline: \emph{General discussion of the NUM problem - Literature of the Back Pressure policy and its dual problem - Difference between the UMW policy and UMW+- Summary of findings in this paper.}
%
Increasingly, networks need to support  a heterogeneous mix of traffic that includes unicast, multicast and broadcast traffic flows.  Such traffic may be comprised of video streaming, file downloads,  distributed computation and storage, as well as a host of other  applications.  The increase in both volume, as well as QoS requirements, of these emerging applications can significantly strain network resources, especially in bandwidth limited wireless networks.   This calls for efficient resource allocation schemes for wireless networks with heterogeneous type of traffic demands.  Unfortunately, however, existing resource allocation schemes, both in theory and practice, are limited to dealing with point-to-point unicast traffic.  While there has been some work on resource allocation for multicast and broadcast flows, until very recently there has been no mechanism for allocating resources in networks with arbitrary traffic flows that include concurrent unicast, multicast and broadcast traffic.  In \cite{report_new}, the authors developed the first such algorithm for joint routing and scheduling in wireless networks with general traffic flows.  The algorithm of \cite{report_new}, called \em{Universal Max Weight} is guaranteed to stabilize the network for any arrivals that are \em{within} the network capacity (stability) region.  However, it provides no service guarantees for traffic that is outside the network capacity region.

When traffic arrivals are outside of the network's capacity region, an admission control mechanism is needed to limit the admitted traffic to being within the network's capacity region.  In particular, a resource allocation scheme is needed that can optimize user-level performance across the network.   Network Utility Maximization (NUM) is an approach for allocating resources across the network in a manner that maximizes overall ``utility,'' where  user's level of satisfaction is measured using a concave function of the allocated data rate to that user (i.e., the utility function).  

Starting with the work of Kelly nearly two decades ago, \cite{kelly-shadowprice},  there has been a tremendous amount of work on the NUM problem (see, for example, \cite{lee-shroff-static-downlink}, \cite{low-duality}, \cite{Mungchiang} and references therein).  These approaches use convex optimization and Lagrangian duality to optimize overall network utility in networks with static traffic demand.   In the context of stochastic traffic, the works of \cite{Neely, Lin, Eryilmaz} used Lyapunov optimization techniques to develop utility optimal flow control and resource allocation for wireless networks with unicast traffic demands.   Unfortunately, however, all of these works deal with unicast (point-to-point) traffic demands.  

In this paper, we develop a jointly optimal admission control, routing, and scheduling algorithm for networks with generalized traffic flows.  Our approach builds upon the Universal Max Weight (UMW) algorithm of \cite{report_new} that solves the optimal routing and scheduling problem for networks with general traffic flows, thus we call our algorithm \emph{UMW+}.   To the best of our knowledge, UMW+ is  the first efficient solution to the NUM problem for networks with generalized traffic flows, including, but not limited to concurrent unicast, broadcast, multicast, and anycast traffic.  We derive UMW+ by using the concept of \emph{precedence-relaxed} virtual queueing network \cite{report_new}. Moreover, by formulating the NUM problem as a static convex optimization problem, we show that, under the dual subgradient algorithm, the evolution of the dual variables precisely corresponds to the evolution of the virtual queues in UMW+. This sheds a new light on the relationship between Max-Weight and sub-gradient algorithms \cite{valls2014relationship}.\edit{We make the following key contributions in this paper:
%%....Need to say more here and then outline the paper.
%       
       \begin{itemize}
      \item     We design the first efficient policy for the NUM problem for \emph{arbitrary} types of traffic, including, but not limited to, concurrent unicast, broadcast, multicast, and anycast sessions. 

       \item    The policy was derived by making effective use of the \emph{virtual network} framework, obtained by relaxing the precedence constraints in a network \cite{report_new}. This methodology was also employed earlier in a wireless broadcasting problem \cite{vq_bcast} and a distributed function computation problem \cite{jianan_paper}. 
  
         \item    We formulate a static version of the NUM problem by decomposing the flows \emph{route-wise}. This is to be contrasted with earlier works on this problem (\emph{e.g.}, \cite{low-duality}, \cite{Mungchiang}), which use node-based flow conservation, applicable \emph{only} to the unicast flows. We show that, under the dual subgradient algorithm with unit step size, the evolution of the dual variables \emph{exactly} corresponds to the evolution of the virtual queues under the proposed UMW+ policy. 
         % This result should be contrasted with a similar result known for the Backpressure policy, which is applicable only to the unicast setting and involves a different problem formulation involving flow conservations at nodes.

            \item     We explicitly characterize the dual objective function of the above optimization problem in terms of cost of the shortest route for each source and the Fenchel conjugates of the utility function.

            \item    Finally, we validate our theoretical results through extensive numerical simulations.

    \end{itemize}
    }

%
%
%\bibitem{lee-shroff-static-downlink}
%J.~W. Lee, R.~R. Mazumdar, and N.~B. Shroff.
%\newblock Downlink power allocation for multi-class cdma wireless networks.
%\newblock {\em IEEE Proceedings of INFOCOM}, 2002.
%
%\bibitem{kelly-shadowprice}
%F.P. Kelly, A.Maulloo, and D.~Tan.
%\newblock Rate control for communication networks: Shadow prices, proportional
%  fairness, and stability.
%\newblock {\em Journ. of the Operational Res. Society}, 49, p.237-252, 1998.
%
%\bibitem{low-duality}
%S.~H. Low.
%\newblock A duality model of tcp and queue management algorithms.
%\newblock {\em IEEE Trans. on Networking}, Vol. 11(4), August 2003.
%
%\bibitem{Mungchiang}
%M Chiang, S. Low, R. Calderbank, J. Doyle
%\newblockLayering as optimization decomposition: A mathematical theory of network architectures
%\newblock {\emProceedings of the IEEE}, Vol. 95 (1), 255-312, 2007.
%
%\bibitem{Neely}
%Fairness and optimal stochastic control for heterogeneous networks
%MJ Neely, E Modiano, CP Li - IEEE/ACM Transactions on Networking (TON), 2008
%
%\binitem{Lin}
%A tutorial on cross-layer optimization in wireless networks
%X Lin, NB Shroff, R Srikant - IEEE Journal on Selected areas in Communications, 2006.
%
%\bibitem{Eryilmaz}
%Fair resource allocation in wireless networks using queue-length-based scheduling and congestion control
%A Eryilmaz, R Srikant - IEEE/ACM Transactions on Networking (TON), 2007
%
%
%

The rest of the paper is organized as follows: In Section \ref{system_model}, we describe the network and traffic model. In Section \ref{virtual_queues}, we detail the virtual network framework and derive the UMW+ control policy. In Section \ref{physical_net}, we prove the stability of the physical queues under the UMW+ policy. In Section \ref{dual_section}, we establish the equivalence between the UMW+ policy and the dual subgradient algorithm. In Section \ref{simulation_section}, we provide simulation results and conclude the paper in Section \ref{conclusion_section}. 

\section{System Model and Problem Formulation} \label{system_model}
\subsection{System Model}
%In this subsection, we formally describe the network model, interference model and the traffic model. 
\paragraph*{Network Model}
 Let the topology of a wireless network be given by the graph $\mathcal{G}(V,E)$ - where $V$ is the set of $n$ nodes and $E$ is the set of $m$ links. Time is slotted, and at every slot, each wireless link could be in either the $\textrm{ON (1)}$ or, $\textrm{OFF (0)}$ state. The random link state process $\bm{\sigma}(t)\in \{0,1\}^{m}$ is assumed to be evolving according to a stationary ergodic stochastic process. If at a slot $t$, link $e$ is \textrm{ON} and is \emph{activated}, it can transmit $c_e$ packets in that slot \footnote{A link does not transmit any packet if it is OFF.}. 
\paragraph*{Interference Model} Due to the wireless inter-channel interference, not all links can be activated at a slot simultaneously. The set of  all interference-free links, which may be activated together in a slot, is given by $\mathcal{M} \subseteq 2^E$. As an example, in the case of \emph{primary} interference constraints, the set $\mathcal{M}$ consists of the set of all \emph{Matchings} in the graph. On the other hand, in an interference-free wired network, the set $\mathcal{M}$ consists of the set of all subsets of links. 
\paragraph*{Traffic and Utility Model} External packets from different traffic \emph{classes} are admitted to the network by an admission controller $\mathcal{A}$. The set of all possible classes of traffic is denoted by $\mathcal{K}$. A traffic class $k \in \mathcal{K}$ is associated with the following two attributes - (1) type of the traffic (\emph{e.g.}, unicast, broadcast, multicast or anycast), and (2) a monotone increasing strictly concave Utility function $U_{(k)}: \mathbb{R}_+\to \mathbb{R}_+$. Let $\mathcal{T}^k$ denote the set of all possible routes for routing packets from class $k$ in the graph $\mathcal{G}$.  As an example, for a unicast class $k$ with source node $s$ and destination node $t$, the set $\mathcal{T}^k$ consists of the set of all $s-t$ paths in the graph. Similarly, for a broadcast class, the set of all routes is given by the set of all spanning trees in the graph. 
%An ON edge, if activated, can serve packets at the rate of $c_e$ packets per slot. 
%and packets of length $s^{(k)}$ bits. Note that, we allow the possibility that different classes of packets may have different lengths, and hence, different service times. 
%In particular, the edge $e$ can serve $c_e^{(k)}= \lfloor\frac{r_e}{s^{(k)}} \rfloor$ packets from class $k$ per slot. 
We consider an \emph{infinitely backlogged} traffic model, where the admission controller $\mathcal{A}$ has potentially an unlimited number of packets from each class available for admission.  
\subsection{Admissible Policies} \label{admissible}
An admissible policy $\pi$ consists of the following three modules: (1) an admission controller $\mathcal{A}$, (2) a routing module $\mathcal{R}$, and (3) a link scheduler module $\mathcal{S}$. The admission controller determines the number of external packets to be admitted to the network from each class in each slot. We assume that, due to physical constraints (power, capacity limitations etc.), at most $A_{\max}$ number of packets may be feasibly admitted from any class per slot \footnote{Taking $A_{\max}\equiv \sum_e c_e$ yields the same optimal utility obtained without this constraint.}. The routing module routes the admitted packets according to the type of flow they belong to, and the link scheduler activates a subset of interference-free links from the set $\mathcal{M}$ in every slot. The set of all admissible policies is denoted by $\Pi$. 
\begin{figure*}[!hbp]
% ensure that we have normalsize text
\normalsize
% Store the current equation number.
%\setcounter{MYtempeqncnt}{\value{equation}}
% Set the equation number to one less than the one
% desired for the first equation here.
% The value here will have to changed if equations
% are added or removed prior to the place these
% equations are referenced in the main text.
%\setcounter{equation}{8}
\hrulefill
%\vspace{4pt}
\begin{eqnarray*}
 \textsc{\textbf{Obj.}}=\underbrace{2\mathbb{E}\bigg(\sum_k A^k(t)\big(\sum_{e} \tilde{Q}_e(t)\mathds{1}(e \in T^k(t))\big)|\bm{\tilde{Q}}(t), \bm{\sigma}(t)\bigg)}_{(a)}
 -\underbrace{2\mathbb{E}\bigg(\sum_e \tilde{Q}_e(t)\mu_e(t) c_e\sigma_e(t)|\bm{\tilde{Q}}(t), \bm{\sigma}(t)\bigg)}_{(b)}
 -\underbrace{2V\sum_k U_k(A^k(t))}_{(c)}, \hspace{4pt} \bm{(*)}
 \end{eqnarray*}
 \hrulefill
% Restore the current equation number.
%\setcounter{equation}{\value{4}}
% IEEE uses as a separator

% The spacer can be tweaked to stop underfull vboxes.
%\vspace*{4pt}
\end{figure*}

\subsection{The Network Utility Maximization (NUM) Problem}
The Network Utility Maximization problem seeks to find an admissible policy $\pi \in \Pi$, which maximizes the sum utility of all classes, while keeping the queues in the network \emph{stable} \footnote{Throughout the paper, by stability, we will mean almost sure rate stability, defined in Eqn. \eqref{stability}. }. Formally, let the random variable $R^\pi_{k}(T)$ denote the number of packets received in common by the destination(s) of class $k$ up to time $T$ under the action of an admissible policy $\pi$. Also, denote the random queue length of packets waiting to cross the edge $e$ under the policy $\pi$ by $Q^\pi_e(T)$. Then the NUM problem seeks to find a policy $\pi^* \in \Pi$, which solves  the following problem:

\begin{equation} \label{NUM_prob}
\max_{\pi \in \Pi} \mathbb{E}\sum_k U_k(r_k)
\end{equation}
Subject to,
\begin{eqnarray}
 &&\lim_{T \to \infty} \frac{R^\pi_{k}(T)}{T} = r_k, \hspace{5pt} \forall k \in \mathcal{K}, \textrm{ w.p. 1.},\\
 &&\lim_{T \to \infty} \frac{1}{T}\sum_{e\in E }Q^{\pi}_e(T)=0, \hspace{5pt} \textrm{ w.p. 1.}, \label{stability}
\end{eqnarray}
where the expectation and the almost sure limits are taken over the randomness of the network configurations $\bm{\sigma}(t)$ and possible randomness in the policy. In the following, we will drop the superscript $\pi$ from the random variables when the driving policy $\pi$ is clear from the context.

\section{The Virtual Network Framework} \label{virtual_queues}
In this section, we describe a \emph{virtual network} framework, obtained by relaxing the natural precedence constraints associated with a multi-hop network. Our goal in this section is to design a utility-optimal stabilizing control policy for the simpler \emph{virtual queueing system}. Section \ref{physical_net} shows that when the same policy is used in the actual physical network, the physical queues are also stable. The virtual network methodology in a similar context was first introduced in \cite{report_new}.\\
%and then employ the same control for the physical queueing system. We will finally show that the resulting hybrid control policy is utility optimal and stabilizing. \\
Consider an admissible policy $\pi$, which admits $A^k(t)$ packets from class $k$, and activates the links $\bm{\mu}(t) \in \{0,1\}^m$ at slot $t$. Taking the random wireless link states into account, the service rate for the virtual queue $\tilde{Q}_e$ is $\mu_e(t)c_e\sigma_e(t)$ packets per slot. Since the virtual network is precedence-relaxed, all admitted packets \emph{immediately} enter all virtual queues on the selected route, i.e., unlike the physical network, arrivals to a virtual queue need not wait to cross the intermediate links \cite{report_new}. Thus, the virtual queues $\tilde{\bm{Q}}(t)$ evolve as %Thus, the dynamics of the virtual queues $\tilde{\bm{Q}}(t)$ is given by
\begin{eqnarray} \label{VQ_dyn}
 \tilde{Q}_e(t+1)=\big(\tilde{Q}_e(t)+\tilde{A}_e(t)-\mu_e(t) c_e\sigma_e(t)\big)^+ \hspace{10pt} \forall e \in E,
\end{eqnarray}

where $\tilde{A}_e(t)$ is the total number of (controlled) arrival of packets to the virtual queue $\tilde{Q}_e$ at slot $t$. Clearly, $\tilde{A}_e(t)$ depends on the routes selected by the routing module $\mathcal{R}$, i.e.,  
\begin{equation}
 \tilde{A}_e(t)=\sum_{k} A^k(t)\mathds{1}(e \in T^k(t)).
\end{equation}
In order to design a utility-optimal stabilizing control policy for the virtual queues, we use the \emph{drift-plus-penalty} framework of \cite{neely2010stochastic}. Consider the following Lyapunov function, which is quadratic in the virtual queue lengths (as opposed to the usual physical queue lengths \cite{neely2010stochastic})
\begin{equation}
 L(\bm{\tilde{Q}}(t))= \sum_e \tilde{Q}_e^2(t).
\end{equation}
The one-slot conditional drift of $L(\bm{\tilde{Q}}(t))$ under the action of a control policy $\pi$ is given as follows
%\hspace*{-10pt}
\begin{eqnarray}
 &&\Delta(\bm{\tilde{Q}}(t), \bm{\sigma}(t))\equiv \mathbb{E}\bigg(L(\bm{\tilde{Q}}(t+1)) - L(\bm{\tilde{Q}}(t))|\bm{\tilde{Q}}(t), \bm{\sigma}(t)\bigg) \nonumber \\
 &&\stackrel{(a)}{\leq}  \mathbb{E}\bigg(\sum_e \big( \tilde{A}_e^2(t) + c_e^2 \nonumber\\
 && +2\tilde{Q}_e(t)(\tilde{A}_e(t)- \mu_e(t) c_e\sigma_e(t))  \big)|\bm{\tilde{Q}}(t), \bm{\sigma}(t)\bigg) \nonumber\\
 &&= B + 2\mathbb{E}\bigg(\sum_e \tilde{Q}_e(t)\sum_{k} \big(A^k(t)\mathds{1}(e \in T^k(t))\big)|\bm{\tilde{Q}}(t), \nonumber\\
 && \bm{\sigma}(t)\bigg)- 2\mathbb{E}\bigg(\sum_e \tilde{Q}_e(t)\mu_e(t) c_e\sigma_e(t)|\bm{\tilde{Q}}(t), \bm{\sigma}(t)\bigg)\nonumber\\
 &&\stackrel{(b)}{=} B+ 2\mathbb{E}\bigg(\sum_k A^k(t)\big(\sum_{e} \tilde{Q}_e(t)\mathds{1}(e \in T^k(t))\big)|\bm{\tilde{Q}}(t), \nonumber\\
 &&\bm{\sigma}(t)\bigg) - 2\mathbb{E}\bigg(\sum_e \tilde{Q}_e(t)\mu_e(t) c_e\sigma_e(t)|\bm{\tilde{Q}}(t), \bm{\sigma}(t)\bigg), \label{drift_expr}
 \end{eqnarray}
where the inequality (a) is obtained by using the virtual queue dynamics in Eqn. \eqref{VQ_dyn}, the equality (b) is obtained by interchanging the order of summation in the first term, $B$ is a finite constant, upper bounded by $kmA_{\max}^2+mc_{\max}^2 $, where $A_{\max}$ is the maximum number of external admissions per slot per class (defined in Section \ref{admissible}) and $c_{\max}\stackrel{\textrm{def}}{=}\max_{e\in E} c_e$.\\
Moreover, admission of $A^k(t)$ packets from class $k$ in slot $t$ yields a ``one-slot utility" of $U_k(A^k(t))$ for class $k$. Following the \emph{``Drift-Plus-Penalty''} framework of \cite{neely2010stochastic}, we consider a cross-layer admission control, routing and link scheduling policy, which is obtained by minimizing the objective function $(*)$ given at the bottom of this page, over all admissible controls $(\bm{A}(t), \bm{T}(t), \bm{\mu}(t))$ per slot.
%\begin{eqnarray} \label{dpp}
% &&\textsc{Per-slot Objective}= \nonumber\\
% &&\underbrace{2\mathbb{E}\bigg(\sum_k A^k(t)\big(\sum_{e} \tilde{Q}_e(t)\mathds{1}(e \in T^k(t))\big)|\bm{\tilde{Q}}(t), \bm{\sigma}(t)\bigg)}_{(a)}\nonumber \\
% &&- \underbrace{2\mathbb{E}\bigg(\sum_e \tilde{Q}_e(t)\mu_e(t) c_e\sigma_e(t)|\bm{\tilde{Q}}(t), \bm{\sigma}(t)\bigg)}_{(b)}\\
% &&-\underbrace{V\sum_k U_k(A^k(t))}_{(c)},\nonumber
%\end{eqnarray}
%\begin{figure*}[!hbp]
%% ensure that we have normalsize text
%\normalsize
%% Store the current equation number.
%%\setcounter{MYtempeqncnt}{\value{equation}}
%% Set the equation number to one less than the one
%% desired for the first equation here.
%% The value here will have to changed if equations
%% are added or removed prior to the place these
%% equations are referenced in the main text.
%\setcounter{equation}{8}
%\begin{eqnarray} \label{dpp}
% \underbrace{2\mathbb{E}\bigg(\sum_k A^k(t)\big(\sum_{e} \tilde{Q}_e(t)\mathds{1}(e \in T^k(t))\big)|\bm{\tilde{Q}}(t), \bm{\sigma}(t)\bigg)}_{(a)}
% -\underbrace{2\mathbb{E}\bigg(\sum_e \tilde{Q}_e(t)\mu_e(t) c_e\sigma_e(t)|\bm{\tilde{Q}}(t), \bm{\sigma}(t)\bigg)}_{(b)}
% -\underbrace{V\sum_k U_k(A^k(t))}_{(c)},\nonumber
%\end{eqnarray}
%% Restore the current equation number.
%%\setcounter{equation}{\value{4}}
%% IEEE uses as a separator
%%\hrulefill
%% The spacer can be tweaked to stop underfull vboxes.
%%\vspace*{4pt}
%\end{figure*}
%\FloatBarrier
For the objective $(*)$, we have substituted the upper bound of the drift $\Delta(\bm{\tilde{Q}}(t), \bm{\sigma}(t))$ from Eqn. \eqref{drift_expr} (without the constant $B$), and $V$ is taken to be a fixed positive constant. This yields the following joint routing, admission control, and link scheduling policy $\pi^*$, which we call \emph{Universal Max-Weight Plus} (\textbf{UMW+}):\\ 
\begin{framed}
\textbf{1. Routing ($\mathcal{R}$):} The routing policy follows by minimizing the term (a) of the objective $(*)$. Consider a weighted graph $\mathcal{\tilde{G}}$, where each edge $e$ is weighted by the  corresponding virtual queue length $\tilde{Q}_e(t)$. Under the policy $\pi^*$, all admitted packets from class $k \in \mathcal{K}$ (refer to $\mathcal{A}$ in Part (2) below) are assigned a route corresponding to the shortest route $T^k(t) \in \mathcal{T}^k$ in $\mathcal{\tilde{G}}$. In particular, 
\begin{itemize}
 \item For a unicast $s-t$ flow, $T^k(t)$ is the weighted shortest $s-t$ path in $\mathcal{\tilde{G}}$.
 \item For a broadcast flow originating from the node $r$, $T^k(t)$ is the Minimum Weight Spanning Tree (MST) in the weighted graph $\mathcal{\tilde{G}}$.
 \item For a multicast flow, $T^k(t)$ is the corresponding Steiner tree in the weighted graph $\mathcal{\tilde{G}}$. 
 \item For an anycast flow \cite{sinha_CDN_journal}, $T^k(t)$ is the weighted shortest path from the source to \emph{any} of the given destinations. 
\end{itemize}
\end{framed}
Denote the cost of the weighted shortest route corresponding to class $k$ obtained above by $C^k(t)$, i.e., 
\begin{equation} \label{shortest_route}
 C^k(t)=\min_{T^k \in \mathcal{T}^k}\sum_e \tilde{Q}_e(t)\mathds{1}(e \in T^k)
\end{equation}

\begin{framed}
 \textbf{2. Admission Control ($\mathcal{A}$):} The admission control policy follows by jointly considering the terms (a) and (c) of the per-slot objective function $(*)$. The \textbf{UMW+} policy admits $A^k(t)$ packets from class $k \in \mathcal{K}$ in slot $t$, where $A^k(t)$ is obtained by the solution to the following one dimensional convex optimization problem:
 \begin{equation} 
  A^k(t) = \argmin_{0\leq x \leq A_{\max}} \big(C_k(t)x -VU_k(x)\big) 
 \end{equation}
\end{framed}
In the case of differentiable utility functions $U_k(\cdot)$, the number of packets $A^k(t)$ admitted may be obtained in closed form as follows: 
\begin{equation}\label{adm_ctrl}
 A^k(t)= \lbrack U_k'^{-1}(C_k(t)/V) \rbrack_{0}^{A_{\max}},
\end{equation}
where $\lbrack \cdot \rbrack_{a}^{b}$ denotes projection onto the interval $\lbrack a,b\rbrack$. \\
\begin{framed}
 \textbf{3. Link Scheduling ($\mathcal{S}$):} The drift-minimizing link scheduling policy is obtained by minimizing the term (b) of the per-slot objective function $(*)$. Consider the weighted graph $\hat{\mathcal{G}}$, where the link $e$ of the graph $\mathcal{G}$ is given a weight of $c_e\tilde{Q}_e(t)$, if it is ON in slot $t$ (i.e. $\sigma_e(t)=1$), and $0$, otherwise. Then the link scheduler $\mathcal{S}$ activates the set of ON links which maximizes the total weight among the set of all interference-free links $\mathcal{M}$, i.e., 
 \begin{equation} \label{max-wt}
  \bm{\mu}^*(t) \in \arg\max_{\bm{\mu}\in \mathcal{M}} \sum_{e} \tilde{Q}_e(t)c_e \mu_e \sigma_e(t).
 \end{equation}

\end{framed}
In the special case of primary interference constraints, where the set $\mathcal{M}$ consists of the set of all \emph{Matchings} of the graph, the problem \eqref{max-wt} corresponds to the Maximum-Weighted Matching problem, which can be solved efficiently even in a distributed fashion \cite{hoepman2004simple}. \\
Our main result in this section is Theorem \ref{VQ_stability}, which claims that under the UMW+ control policy described above, the virtual queues are stable and the average expected utility obtained is arbitrarily close to the optimal utility.

 \begin{framed}
 \begin{theorem} \label{VQ_stability}
  Let $U^*$ be the optimal utility for the NUM problem given in Eqn. \eqref{NUM_prob}. Then, under the action of the \textbf{UMW+} control policy, (a) the virtual queues are rate stable (in the sense of Eqn.\eqref{stability}), and (b) the utility achieved is at least $U^*-\mathcal{O}(\frac{1}{V})$, for any $ V>0$.
 \end{theorem}
\end{framed}
\begin{proof}
 See the Appendix (Section VIII). 
% See Appendix \ref{optimality_proof}.\\
\end{proof}
An alternative treatment of Utility Optimality and Stability of the Virtual Queues under the UMW+ policy will be given in Section \ref{dual_section}, where we relate the virtual queue evolution to the subgradient descents of an appropriately defined dual optimization problem. This sheds new and fundamental insights in the construction and operation of the virtual queues.  

\section{Control of the Physical Network} \label{physical_net}
In the physical network, \emph{the same} admission control, routing and link scheduling policy are used as in the virtual network of Section \ref{virtual_queues}. 
Hence, the same number of packets are admitted at each slot from each class to the virtual and the physical network. As a result, the utility achieved in these two networks are the same, however, their queueing evolutions are different. In this Section, we establish the stability of the physical queues, which, unlike the virtual queues, are subjected to the usual precedence constraint of a multi-hop network. 
\subsection{Packet arrivals to the Physical Queues}\label{phy_Q_arr}
To connect the arrivals of packets in the virtual and the physical network, observe that since the routes are fixed at the sources, the total number of physical packets $A_e(t_1,t_2)$ admitted to the physical network in any time interval $(t_1,t_2]$, that \emph{will cross edge $e$ in future}, is the same as the number of virtual packets arrivals $\tilde{A}_e(t_1,t_2)$ to the corresponding virtual queue $\tilde{Q}_e$.  Let the total service allocated to serve the virtual queue (resp. physical queue) in the time interval $(t_1,t_2]$ be denoted by $\tilde{S}_e(t_1,t_2)$ (resp. $S_e(t_1,t_2)$). Using the Skorokhod map \cite{kleinrock1975queueing} for the virtual queue iterations \eqref{VQ_dyn}, we have 
\begin{eqnarray} \label{skorokhod}
	\tilde{Q}_e(t)= \sup_{0\leq t_1 \leq t}\bigg(\tilde{A}_e(t_1,t)-\tilde{S}_e(t_1,t)\bigg).
\end{eqnarray}
However, as noted above, we have
\begin{eqnarray} \label{connection}
	A_e(t_1,t_2)=\tilde{A}_e(t_1,t_2), S_e(t_1,t_2)= \tilde{S}_e(t_1,t_2).
\end{eqnarray}
Hence, using Theorem \eqref{VQ_stability} for the stability of the virtual queues and combining it with Eqns \eqref{skorokhod} and \eqref{connection}, we have 
\begin{eqnarray} \label{arr_bound}
	A_e(t_1,t) \leq S_e(t_1,t) + M(t), \forall t_1<t, \emph{a.s.},
\end{eqnarray}  
where $M(t)= o(t)$. 
\subsection{Stability of The Physical Queues} \label{phy_Q_stability}
Note that, Theorem \ref{VQ_stability} establishes the stability of the \emph{virtual queues} under the action of the \textbf{UMW+} policy. To prove the optimality of this policy in the NUM setting of Eqn. \eqref{NUM_prob}, we need to establish
%relate the stability of the virtual queues from Theorem \ref{VQ_stability} with the NUM problem, which is concerned with the physical network with the usual precedence constraints. In particular, this requires us to 
%connect the stability of the virtual queues $\{\tilde{\bm{Q}}(t)\}_{t \geq 1}$ to the 
 the stability of the physical queues $\{\bm{Q}(t)\}_{t \geq 1}$. This is a non-trivial task because, unlike the virtual queues, the dynamics of the physical queues is subject to the precedence constraints, and regulated by the \emph{packet scheduling policy} employed in the physical network.  A packet scheduling policy is a rule which resolves the contention when multiple packets want to cross the same edge at the same slot. 
%\end{framed}
The most common examples of packet scheduling policies include First In First Out (\textbf{FIFO}), Last In First Out (\textbf{LIFO}) disciplines etc.  Following our earlier work \cite{report_new}, we consider a simple packet scheduling policy called the Extended Nearest to Origin (\textbf{ENTO}) policy:
\begin{definition}[Extended Nearest to Origin policy]
The ENTO policy prioritizes packets according to the decreasing order of their current distances (measured in hop-lengths) from their respective sources. 
\end{definition}
For example, if there are two packets $p_1$ and $p_2$, which have traversed, say, $10$ and $20$ hops respectively from their sources, and wish to cross the same active edge $e$ (with a unit capacity per slot) in the same slot, the ENTO policy prioritizes the packet $p_1$ over $p_2$ to cross the edge $e$.\\  
Using ideas from adversarial queueing theory \cite{gamarnik}, it is shown in \cite{report_new} (Theorem 3) that, under the ENTO packet scheduling policy, the almost sure packet arrival bound \eqref{arr_bound} also implies the rate stability of the physical queues. Combining this result with Theorem \ref{VQ_stability}, we conclude the following:
\begin{framed}
\begin{theorem}[Stability of the Physical Queues] \label{Phy_Q_stability}
 Under the action of the \textbf{UMW+} control policy, the physical queues are rate stable, i.e.,
 \begin{eqnarray*}
  \lim_{T \to \infty} \frac{1}{T}\sum_{e\in E }Q_e(T)=0, \hspace{5pt} \textrm{ w.p. 1.}
 \end{eqnarray*}

\end{theorem}
\end{framed}
The proof of Theorem \ref{Phy_Q_stability} follows directly from the proof of Theorem 3 of \cite{report_new} and uses the bound \eqref{arr_bound} as a starting point. 
Combining Theorem \ref{VQ_stability} with Theorem \ref{Phy_Q_stability}, we conclude that the proposed \textbf{UMW+} policy is a utility optimal stabilizing control policy, which solves the NUM problem efficiently. This concludes the first half of the paper.

\section{A Dual Perspective on the Virtual Queue Dynamics} \label{dual_section}
In this Section, we consider the dual of an offline version of the NUM problem in a static wired network setting (no interference constraints). We give an alternative derivation of the UMW+ policy  from an optimization theory perspective, which sheds further insight into the structure of the optimal policy. Our motivation in this section is similar to \cite{Mungchiang}, \cite{Lin}, which give similar development for the Backpressure policy \cite{tassiulas}. The dual problem is also of sufficient theoretical and practical interest, as the powerful machinery of convex optimization may be used to derive alternative efficient algorithms for solving the dual problem, which may then be translated to other dynamic policies (apart from the UMW+) for solving the NUM problem. 

 %To simplify notations, we may assume without any loss of generality (after relabelling, if necessary) that the set of \emph{routes} $\{\mathcal{
%T}_k, k=1,2, \ldots, |\mathcal{K}|\}$ are disjoint. 
As in the previous section, the topological structure of a route depends on the type of flow - \emph{e.g.}, a route is a path for unicast and anycast flows, a spanning tree for broadcast flows, a Steiner tree for multicast flows etc. \\
Fix a strictly positive constant $V$. Our goal is to solve the following utility maximization problem $\mathcal{P}$:
\begin{framed}
\begin{equation} \label{primal}
\hspace{-70pt}\textrm{Problem } \mathcal{P}: \hspace{10pt} \max V\sum_k U_k(r_k)
\end{equation}
Subject to, 
\begin{eqnarray}
 r_k &=& \sum_{p \in \mathcal{T}_k}f_p, \hspace{10pt} \forall k \in \mathcal{K}.\label{decomp}\\
 \sum_{p : e \in p} f_p &\leq& c_e, \hspace{10pt} \forall e \in E.  \label{capacity_constr}\\
 \bm{f}&\geq& \bm{0} \label{nn}.
\end{eqnarray}
\end{framed}
The objective \eqref{primal} denotes the total utility, scaled by $V$. The constraint \eqref{decomp} is obtained by decomposing the total incoming flow $r_k$ to a class $k$ into all available routes, where the variable $f_p$ corresponds to the amount of flow carried by the route $p \in \mathcal{T}^k$. The constraint \eqref{capacity_constr} corresponds to the capacity of edge $e$, and the constraint \eqref{nn} corresponds to the non-negativity property of the flow variables. \\
In the special case of the NUM problem dealing exclusively with the unicast flows, the flow decomposition constraint \eqref{decomp} is usually replaced with the \emph{flow conservation} constraints at the nodes, which leads to the Backpressure policy \cite{Mungchiang}. In contrast, we formulate the NUM problem with the flow decomposition constraint, since flow, in general, is not conserved at the nodes due to packet replications (as in broadcast and multicast flows). 
\subsection{The Dual Problem $\mathcal{P}^*$}
 The problem $\mathcal{P}$ is a concave maximization problem with linear constraints. To obtain its dual problem, we relax the capacity constraints \eqref{capacity_constr} by associating a non-negative dual variable $q_e$ with the constraint corresponding to the edge $e$. We choose \emph{not to relax} the flow decomposition constraints \eqref{decomp} and the non-negativity constraints \eqref{nn}\footnote{Relaxation of these constraints yields a different dual problem.}. This yields the following partial Lagrangian: 
\begin{equation}\label{lagrangian}
 L(\bm{r}, \bm{f}, \bm{q})= V\sum_k U_k(r_k) + \sum_e q_e (c_e - \sum_{p: e \in p }f_p),
\end{equation}
leading to the following dual objective function $D(\bm{q})$:
\begin{framed}
\begin{equation} \label{dual_obj}
 D(\bm{q}):= \max_{\bm{r}, \bm{f \geq 0}} L(\bm{r}, \bm{f}, \bm{q}),
\end{equation}
Subject to,
\begin{equation}
  r_k = \sum_{p \in \mathcal{T}_k}f_p, \hspace{10pt} \forall k \in \mathcal{K}\label{decomp2}.
\end{equation}
\end{framed}
By strong duality \cite{bertsekas_convex}, the problem $\mathcal{P}$ is equivalent to the following dual problem:

\begin{equation} \label{dual_prob1}
 \hspace{-50pt}\textrm{Problem } \mathcal{P}^*:\hspace{10pt}\min_{\bm{q\geq 0}}D(\bm{q}).
\end{equation}
We next establish a simple lemma which will be useful for our subsequent development.
\begin{framed}
\begin{lemma} \label{sp}
For any fixed $\bm{q \geq 0}, \bm{r} \geq \bm{0}$, an optimal solution to the problem \eqref{dual_obj} is obtained by routing the entire flow $r_k$ from each class $k$ along a weighted shortest route $p_k^*\in \mathcal{T}^k$, weighted by the corresponding dual variables $\bm{q}$. 
\end{lemma}
\end{framed}
%\begin{proof}
\textbf{Proof:}
Exchanging the order of summation in the last term of the Lagrangian in Eqn. \eqref{lagrangian}, we have 
\begin{equation}  \label{dual_new}
 L(\bm{r,f,q})= V\sum_k U_k(r_k) + \sum_e q_e c_e - \sum_{p } f_p \big(\sum_{e:e \in p} q_e \big).
\end{equation}
Define $c_p(\bm{q})\equiv \sum_{e:e \in p} q_e$ to be the cost of the route $p$ where each edge $e$ is weighted by the dual variable $q_e$, $\forall e \in E$. Hence, from the last term of the above expression \eqref{dual_new} and the flow decomposition constraint \eqref{decomp2}, it immediately follows that the objective \eqref{dual_obj} is maximized by routing the entire incoming flow $r_k$ along a weighted shortest route $p_k^* \in \mathcal{T}_k$ corresponding to class $k$.  In other words, for any class $k$, 
%we have $f^*_{p_k^*}=r_k, f^*_p=0, \forall p \in \mathcal{T}^k \neq p_k^*$, 
\begin{eqnarray*}
	f^*_p=\begin{cases}
		r_k, \textrm{ if } p=p_k^*\\
		0, \textrm{ if } p \in \mathcal{T}^k \textrm{ and } p\neq p_k^*,
	\end{cases}
\end{eqnarray*}
where $p_k^*=\arg\min_{p \in \mathcal{T}^k}c_p(\bm{q})$ (ties are broken arbitrarily). $\blacksquare$ \\\\
% it is always optimal to set $f_p=0$ for all paths apart from the shortest route corresponding to the class $k$. This can be seen by feasibly diverting flows from non-shortest routes to a single shortest route. 
%\end{proof}
Define $c_k^*(\bm{q})$ to be the cost of the shortest route corresponding to class $k$, i.e., $ c_k^*(\bm{q})=\min_{p \in \mathcal{T}^k} c_p(\bm{q}).$
%\begin{equation} \label{vq_sp}
% c_k^*(\bm{q})=\min_{p \in \mathcal{T}^k} c_p(\bm{q}). 
%\end{equation}
Since, $c_k^*(\bm{q})$ is defined to be the point wise minimum of several linear functions, it is a concave function of $\bm{q}$ \cite{bertsekas_convex}. For several important traffic classes (\emph{e.g.,} unicast, broadcast, anycast) there are standard combinatorial algorithms for efficiently computing $c_k^*(\bm{q})$ (\emph{e.g.} Weighted Shortest Path, Minimum Weight Spanning Tree etc.). \\
With an optimal setting of the flow variables $\bm{f}$ resolved by Lemma \eqref{sp}, Eqn. \eqref{dual_new} implies that the computation of the dual objective function \eqref{dual_obj} reduces to optimizing the traffic admission rates $r_k$ for each class $k$ as follows: 
\begin{equation}  \label{rate_fn}
r_k^*(\bm{q})= \arg \max_{r_k \geq 0} \big(VU_k(r_k)-r_kc_k^*(\bm{q})\big)
\end{equation}

Due to strict concavity of the utility functions $U_k(\cdot)$, the optimal solution to the problem \eqref{dual_obj} is obtained by setting the derivative of the objective with respect to the variable $r_k$ to zero, which yields
\begin{equation} \label{opt_rate}
 r_k^*(\bm{q})= \bigg(U_k'^{-1}\big(c_k^*(\bm{q})/V\big)\bigg)^+,
\end{equation}
 where we project $r_k$ on the set of non-negative real numbers due to non-negativity constraints of the rates. Substituting Eqn. \eqref{opt_rate} into Eqn. \eqref{dual_obj}, we obtain an implicit expression of the dual objective function $D(\bm{q})$. In the following, we derive an explicit expression of the dual function in terms of the Fenchel conjugate of the utility functions \cite{bertsekas_convex}.
 \subsection{Derivation of the Dual Objective Function} \label{dual_derivation}
%In this section, we derive an explicit expression for the dual objective function $D(\bm{q})$, given in Eqn. \eqref{dual_obj}. 
Substituting the value of $r_k^*(\bm{q})$ into the Lagrangian \eqref{lagrangian} and noting that $f_p^*=r_k^*(\bm{q})$ only along the corresponding shortest route and is zero otherwise, we have 
\begin{equation} \label{dual_fun}
 D(\bm{q})= V\sum_k \big(U_k(r_k^*(\bm{q}))-r_k^*(\bm{q}) \frac{c_k^*(\bm{q})}{V}\big) + \sum_e q_e c_e .
\end{equation}
%(Verify that $D(\bm{q})$ is a convex function of $\bm{q}$). \\
The above expression may also be written in terms of the Fenchel's conjugate \cite{bertsekas_convex} of the utility functions. For this, we recall the definition of the Fenchel conjugate $f^\dagger$ of (an extended real-valued) function $f$: 
\begin{equation}
 f^\dagger(z)= \sup_{x \in \textrm{dom}(f)}\big( x^Tz - f(x) \big)
\end{equation}
Now, let the function $U_k^\dagger:\mathbb{R}\to \mathbb{R}$ denote the Fenchel conjugate of the function $-U_k(\cdot)$, which, by our assumption, is a strictly convex function. Thus, from Eqn. \eqref{dual_fun} we can write
\begin{framed}
\begin{equation} \label{dq_gen}
 D(\bm{q})= V\sum_k U_k^\dagger\big(-\frac{c_k^*(\bm{q})}{V}\big) + \sum_e q_ec_e.
\end{equation}
\end{framed}
%From the above Eqn., we can readily verify that $D(\bm{q})$ is a convex function of $\bm{q}$. 
In the following, we use Eqn. \eqref{dq_gen} to derive explicit functional forms of the dual objective functions for two important classes of network utility functions.\\ 
%\subsubsection{Examples}\\
\textbf{1. Logarithmic Utility Functions:}
Consider the class of Logarithmic utility functions defined as follows:
\begin{eqnarray}
 U_k(r_k)= \gamma_k\log(1+r_k), \hspace{10pt} r_k\geq 0,
\end{eqnarray}
where $\gamma_k$ is a positive constant. Among its many attractive properties, the Logarithmic utility functions ensures proportionally fair rate allocations among all participating classes.  \\
The Fenchel's conjugate of $-U_k(r_k)$ may be computed as 
\begin{eqnarray*}
 U_k^\dagger(z)&=& \sup_{x \geq 0}\big(xz+ \gamma_k \log(1+x)\big)\\
 &=& \begin{cases}
      \gamma_k \log(-\frac{\gamma_k}{z}) - (\gamma_k+z), \textrm{ if } z<0\\
      +\infty, \textrm{ if } z \geq 0
     \end{cases}
\end{eqnarray*}
Hence, the dual function is given as $D(\bm{q})=$
\begin{eqnarray*} 
\hspace{-10pt}\begin{cases}
             V\sum_k \big(\gamma_k \log(\frac{\gamma_kV}{c_k^*(\bm{q})})-\gamma_k+ \frac{c_k^*(\bm{q})}{V}\big) + \sum_e q_ec_e, \textrm{ if } c_k^*(\bm{q}) > 0 \\
             +\infty, \textrm{ if } c_k^*(\bm{q}) \leq 0,
            \end{cases}
\end{eqnarray*}
%where, we recall that $c_k^*(\bm{q})$ is the length of the shortest route for class $k$ in the graph $\mathcal{G}(V,E)$, where its edges are weighted by the dual variables $\bm{q}$.\\
\textbf{2. $\alpha$-fair Utility Functions:}
Next, we consider the $\alpha$-fair utility functions defined as follows:
\begin{equation}
  U_k(r_k)= \gamma_k\frac{r_k^{1-\alpha}}{1-\alpha}, \hspace{10pt} r_k\geq 0,
\end{equation}
where $\gamma_k > 0$ and $0< \alpha< 1$ are positive parameters. \\
The Fenchel's conjugate of $-U_k(r_k)$ may be computed as 
\begin{eqnarray*}
 U_k^\dagger(z)&=& \sup_{x \geq 0}\bigg(xz+ \gamma_k \frac{x^{1-\alpha}}{1-\alpha}\bigg)\\
 &=& \begin{cases}
      \frac{\alpha}{1-\alpha} \gamma_k^{1/\alpha} (-z)^{1-1/\alpha}, \textrm{ if } z<0\\
      +\infty, \textrm{ if } z \geq 0
     \end{cases}
\end{eqnarray*}
Hence, the dual function $D(\bm{q})$ is given as 
\begin{eqnarray*}
&&D(\bm{q})=\\ 
&&\begin{cases}
             V\frac{\alpha}{1-\alpha}\sum_k \gamma_k^{1/\alpha}\big(\frac{c_k^*(\bm{q})}{V}\big)^{1-1/\alpha} + \sum_e q_ec_e, \textrm{ if } c_k^*(\bm{q}) > 0 \\
             +\infty, \textrm{ if } c_k^*(\bm{q}) \leq 0,
            \end{cases}
\end{eqnarray*}
%where, we recall that $c_k^*(\bm{q})$ is the length of the shortest route for class $k$ in the graph $\mathcal{G}(V,E)$, where its edges are weighted by the vector $\bm{q}$.
 \subsection{Subgradient Method and its Equivalence with UMW+}
 Since the dual objective $D(\bm{q})$, as given in Eqn. \eqref{dual_obj}, is a point wise maximum of linear functions, it is convex \cite{bertsekas_convex}. Moreover, the objective $D(\bm{q})$, as seen from Eqn. \eqref{dq_gen}, is \emph{non-differentiable}, as the shortest path cost $c_k^*(\bm{q})$ is not a differentiable function of $\bm{q}$, in general. Hence, we  use a first order method suitable for non-smooth objectives, known as the \emph{Subgradient Descent}  \cite{bertsekas_convex}, to solve the dual problem \eqref{dual_prob1}.\\ 
 A subgradient $\bm{g} \in \partial D(\bm{q})$ of the dual function $D(\bm{q})$ may be computed directly from the capacity constraint of the primal problem (Eqn. \eqref{capacity_constr}) as follows: 
 \begin{equation} \label{subgrad}
  g_e(\bm{q})= c_e -A_e(\bm{q}), \hspace{10pt} \forall e \in E,
 \end{equation}
where $A_e(\bm{q})\equiv \sum_{p: e \in p}f_p^*(\bm{q})=\sum_{k} r_k^*(\bm{q}) \mathds{1}(e \in p_k^*(\bm{q}))$.\\
% obtained from the capacity constraint of the primal problem (Eqn. \eqref{capacity_constr}).\\
 Finally, we solve the dual problem $\mathcal{P}^*$ by the dual subgradient method with constant step-size equal to $\theta>0$. This yields the following iteration for the dual variables $\bm{q}$:
\begin{eqnarray} \label{dual_update}
 q_e(t+1)=(q_e(t) + \theta(A_e(\bm{q})-c_e))^+, \hspace{10pt} \forall e \in E.
\end{eqnarray}

For the step-size $\theta=1$, the above iteration \emph{corresponds exactly} to the virtual queue dynamics. More precisely, Eqn. \eqref{shortest_route} corresponds to Lemma \ref{sp}, which is concerned with routing along the shortest route in the weighted graph; Eqn. \eqref{adm_ctrl} corresponds to Eqn. \eqref{opt_rate}, which corresponds to the packet admission control, and the virtual queue dynamics in Eqn. \eqref{VQ_dyn} corresponds to the dual subgradient update of Eqn. \eqref{dual_update}. Moreover, for any constant step-size $\theta>0$, at any step $t$, the virtual queues $\bm{\tilde{Q}}(t)$ under UMW+ are  exactly equal to the $1/\theta$-scaled versions of the corresponding dual variables $\bm{q}(t)$ under the subgradient descent iterations.
%With the additional assumption that the maximization in Eqn. \eqref{rate_fn} is performed in the bounded range of the variable $[0, A_{\max}]$, for some suitably large positive constant $A_{\max}$, \footnote{This assumption parallels our assumption of bounded arrival per slot in Section \ref{admissible}} we have the following classical convergence result of the subgradient iterations \eqref{dual_update} with constant step-size $\theta>0$:
\paragraph*{Convergence of the Subgradient Descent}
To establish the convergence of the subgradient iterations, it is useful to add an additional constraint in the primal problem $\mathcal{P}$, without changing the optimal solution. For this, observe that, for any feasible solution to $\mathcal{P}$, from the capacity constraint \eqref{capacity_constr}, we have 
\begin{eqnarray*}
	\sum_k r_k \leq \sum_e c_e.
\end{eqnarray*}
With this additional constraint, we conveniently restrict the range of each of the admission control variables $r_k$, in the sub-problems  \eqref{rate_fn}, to $[0, \sum_e c_e]$. This upper-bounds the rate of flows crossing the edge $e$ to $A_e(\bm{q}) \leq |\mathcal{K}|\sum_e c_e \stackrel{\textrm{def}}{=} \bar{A}_{\max}$. Hence, the norms of the subgradients $\bm{g}(t)$ (from Eqn. \eqref{subgrad}) are \emph{uniformly} bounded by 
\begin{eqnarray*}
	||\bm{g}(t)||_2^2 \leq \sum_e (c_e^2 + A^2_e(\bm{q}(t))) \leq m(c_{\max}^2+ m\bar{A}_{\max}^2).
\end{eqnarray*}
We have the following convergence result:
\begin{framed}
\begin{theorem}[Convergence within a Neighborhood] \label{subgrad_convergence}
Under the subgradient iterations \eqref{dual_update}, we have 
\begin{eqnarray*}
VU^* + \frac{\theta b^2}{2} \geq \limsup_{t \to \infty} D(\bm{q}(t)) \geq \liminf_{t \to \infty}D(\bm{q}(t))\geq  VU^*, 
 \end{eqnarray*}
 where $b^2\equiv m(c_{\max}^2+ m|\mathcal{K}|^2(\sum_e c_e)^2)$.
\end{theorem}
\end{framed}

Theorem \ref{subgrad_convergence} follows from an application of Proposition 2.2.2 of \cite{bertsekas2015convex}. This shows that the utility achieved by the dual algorithm is within an additive gap of $O(1/V)$ from $U^*$.
%\textcolor{blue}{We have the following theorem on the (\textbf{Quote Bertsekas Subgradient Convergence Result.})}

% 
% \textbf{Note:} Interestingly, the optimal primal solution does not follows from the dual problem as the Lagrangian is not strictly convex in the primal variable $\bm{f}$. This can be seen from the fact that, all paths apart from the shortest path carry zero flow in any iteration. Nevertheless an optimal primal solution may be obtained (?) by averaging the primal iterates. 

%\input{implementation}
\section{Numerical Simulations} \label{simulation_section}
In this Section, we provide simulation results to explore the performance of the UMW+ policy in diverse network and traffic settings.\\
\subsection{Unicast Traffic in a Wired Network}
\begin{figure}
%\label{unicast_fig}
\hspace{15pt}
	\begin{minipage}{0.1\textwidth}
		\begin{overpic}[scale=0.33]{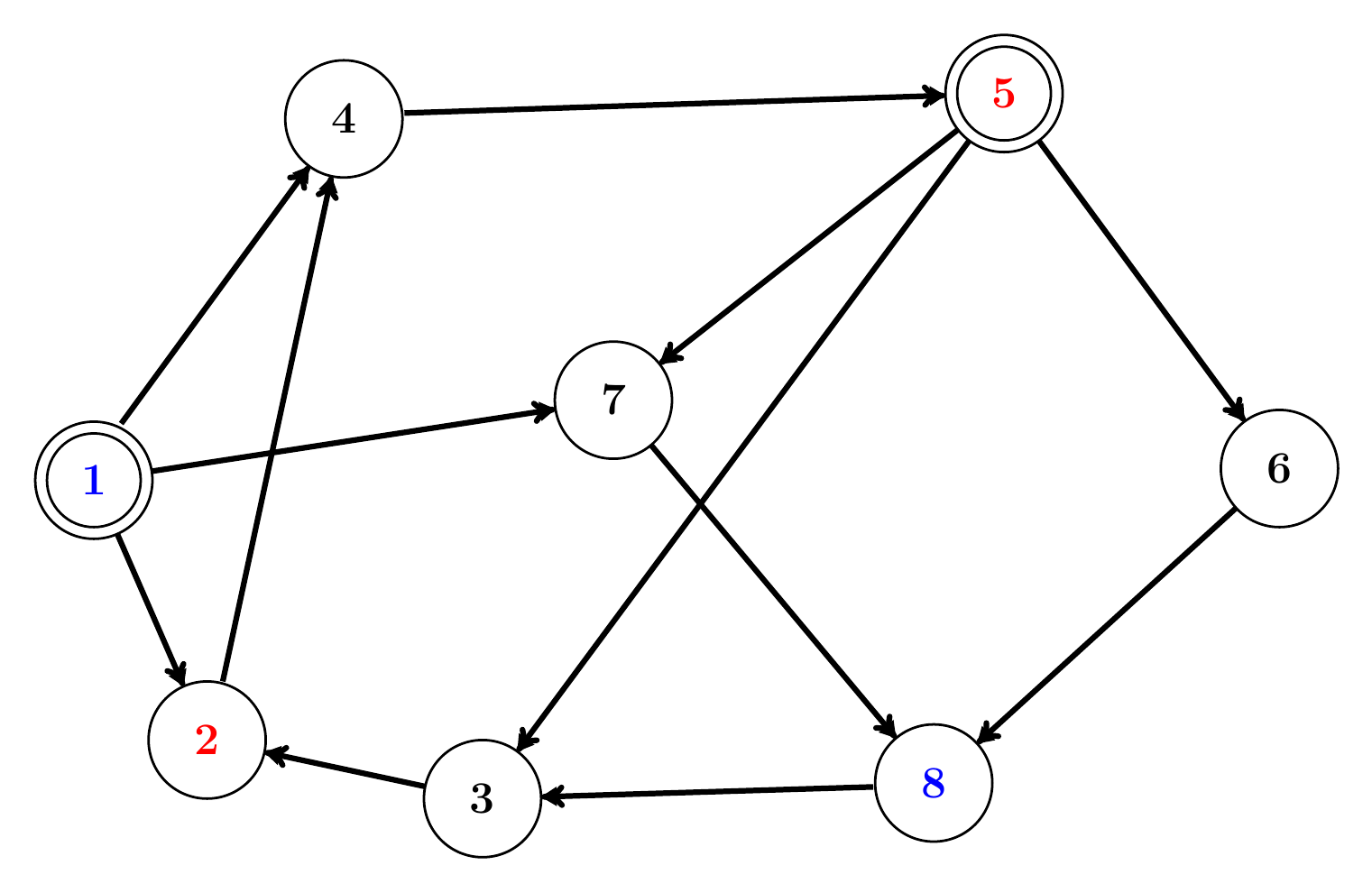}
		\put(-2,35){$\textcolor{blue}{s_1}$}
		\put(-7.6,28){$\textcolor{blue}{\Longrightarrow}$}
		\put(68,-1){$\textcolor{blue}{t_1}$}
		\put(75,65){$\textcolor{red}{s_2}$}
		\put(78,57){$\textcolor{red}{\Longleftarrow}$}
		\put(-13,29){$\textcolor{blue}{r_1}$}
		\put(89,57){$\textcolor{red}{r_2}$}
		\put(15,02){$\textcolor{red}{t_2}$}
		\put(71.5,1){$\textcolor{blue}{\searrow}$}
		\put(7,03){$\textcolor{red}{\swarrow}$}
		\put(45,-4){\small{(a)}}
		\end{overpic}
		%\caption{\small{Wired Network Topology for Unicast Traffic}}
	\end{minipage}
	\hspace{85pt}
	\begin{minipage}{0.14\textwidth}
	 \begin{overpic}[scale=0.4]{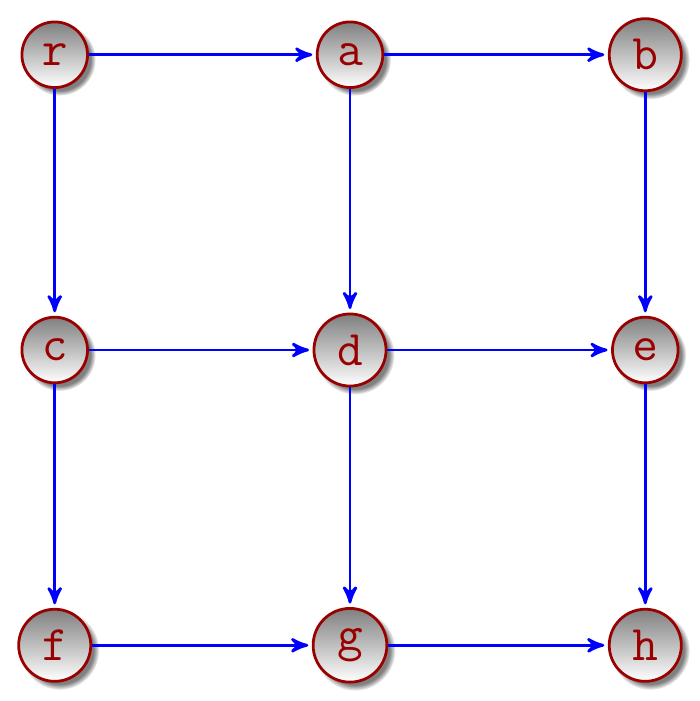}
	 \put(45,-12){\small{(b)}}
	  \end{overpic}
\end{minipage}
\vspace{5pt}
\caption{\small{Network Topologies used for (a) Unicast and (b) Broadcast Simulations}}
\label{Network_Topolgy}
\end{figure}

\begin{figure}
%	\center
	%\includegraphics{0.49\textwidth}
	\hspace{-10pt}
	\begin{minipage}{0.1\textwidth}
		\begin{overpic}[scale=0.26]{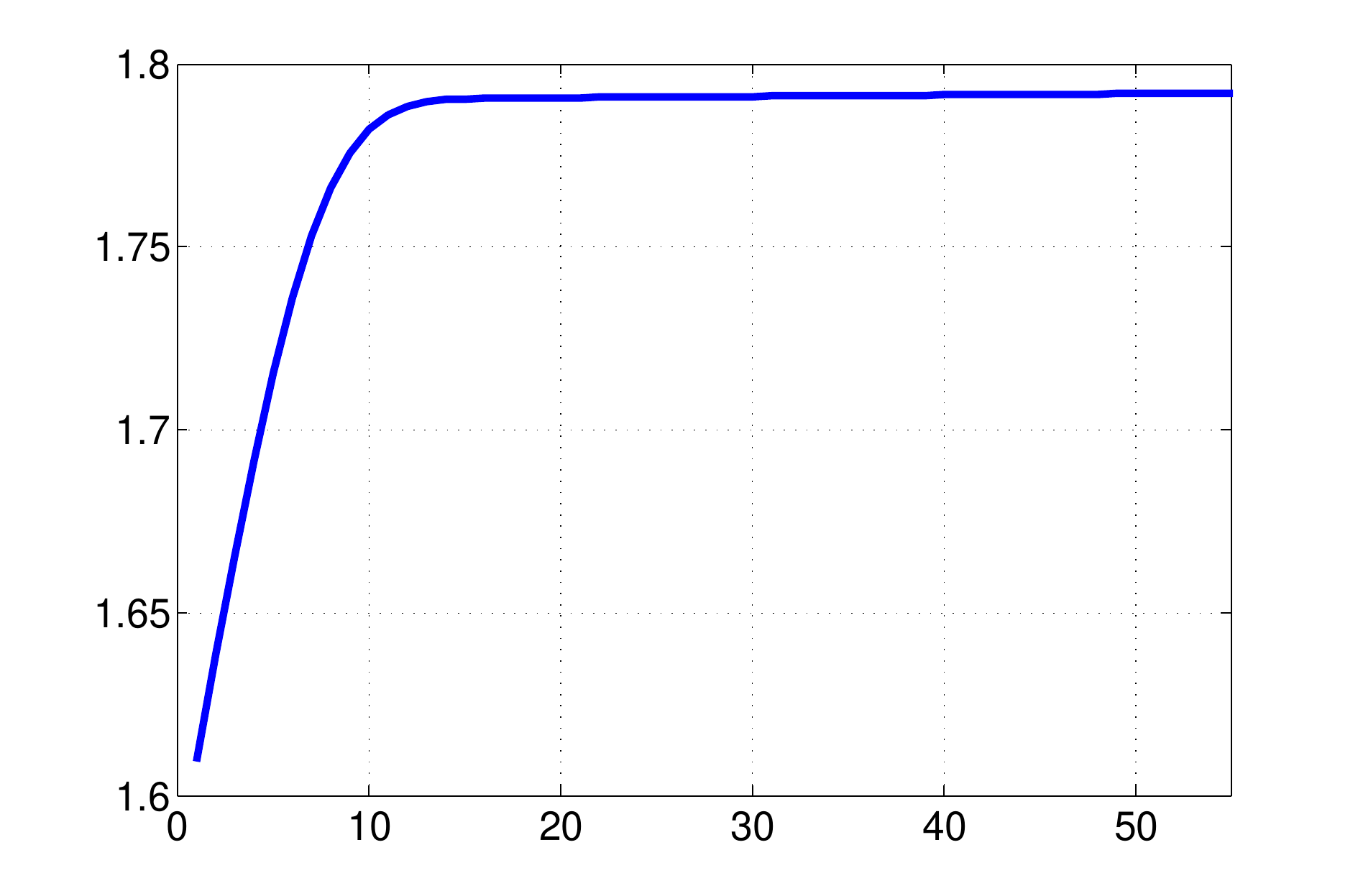}
		%\put(45,40){\small{$U_1(r_1)=\ln(1+r_1)$}}
		\put(2,20){\rotatebox{90}{\scriptsize{Sum Utility}}}
		\put(50,-1){\scriptsize{$V$}}
		\put(48,-7){\small{(a)}}
		\end{overpic}
		%\end{figure}
		%\begin{figure}
		\end{minipage}
		\hspace{70pt}
         \begin{minipage}{0.1\textwidth}
	\begin{overpic}[scale=0.26]{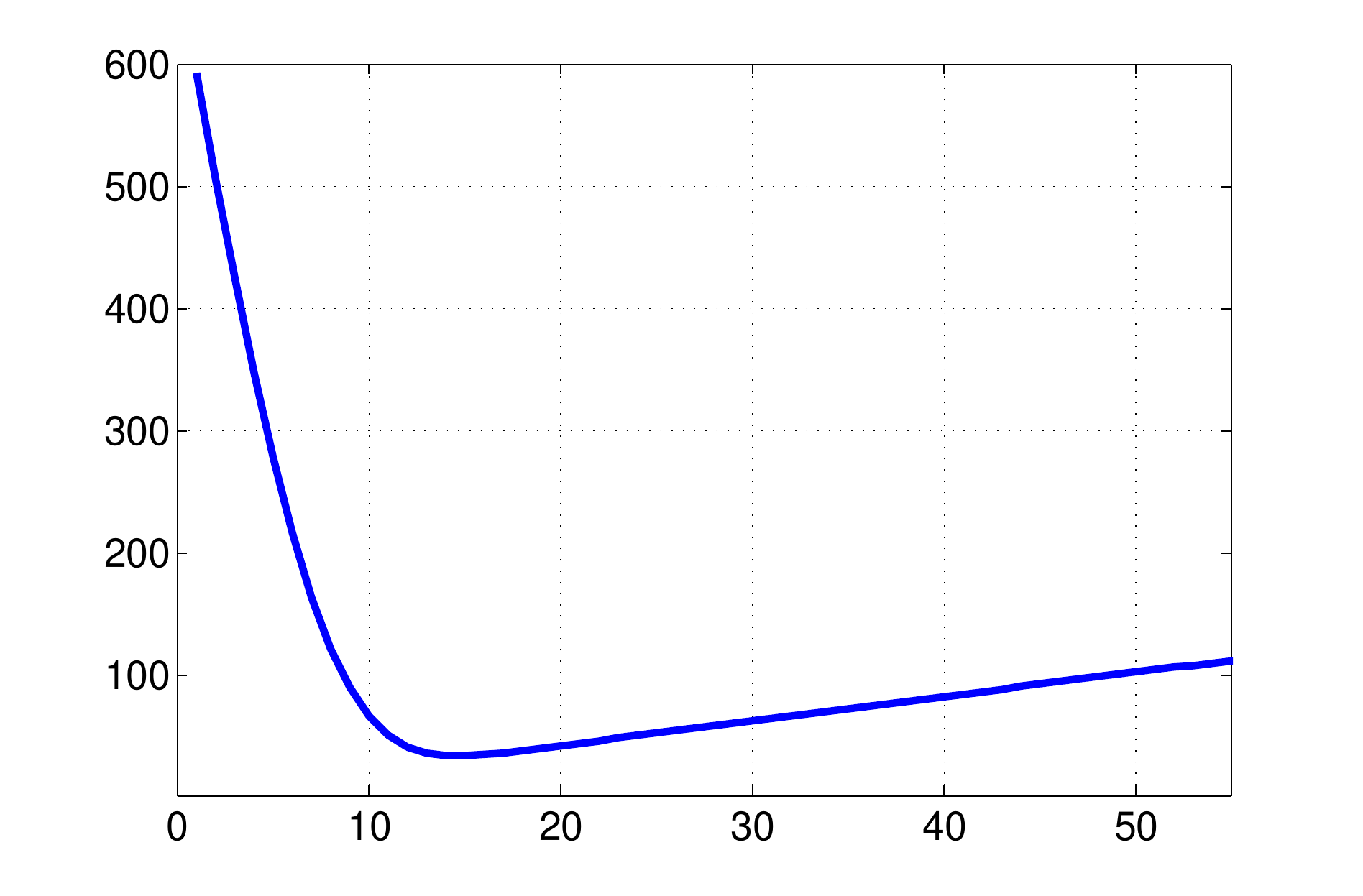}
			\put(2,10){\rotatebox{90}{\scriptsize{Average Queue Lengths}}}
			\put(50,-1){\scriptsize{$V$}}
			\put(48,-7){\small{(b)}}
	\end{overpic}
	\end{minipage}
	\vspace{5pt}	
	\caption{\small{\footnotesize{Performance of the UMW+ policy for the Utility Maximization problem (NUM). }}}
	\label{util_fig}
	%\end{overpic}
\end{figure}

To begin with, we consider the same wired network topology as in \cite{report_new}, with two unicast sources $s_1, s_2$ and destinations $t_1, t_2$, shown in Figure \ref{Network_Topolgy} (a). The capacity of each directed link is taken to be one packet per slot. From the network topology, it can be easily seen that there are two $s_1\to t_1$ paths (\emph{e.g.}, $1\to 4 \to 5 \to 6 \to 8$ and $1 \to 7 \to 8)$ and one $s_2 \to t_2$ path (\emph{e.g.}, $5 \to 3 \to 2$), which are \emph{mutually edge disjoint}. Moreover, we also have $\textsc{cut}(\{ 1,2,3,4,5,6,7\}, \{8\})=2$ and $\textsc{cut}(\{4,5,6,7,8\}, \{1,2\})=1$.  This implies that the optimal solution to the NUM problem is attained at $r_1=2, r_2=1$ for any non-decreasing utility function. \\
Next, we run the proposed UMW+ policy with logarithmic utility functions $U_1(r_1)=\ln(1+r_1), U_2(r_2)=\ln(1+r_2)$. The theoretical optimal utility value, in this case, is easily computed to be $U^*= \ln(3)+\ln(2)\approx 1.79.$ Figure \ref{util_fig} (a) shows the achieved utility by the UMW+ policy with the variation of the $V$ parameter. Figure \ref{util_fig} (b) shows the variation of the total queue length as a function of $V$. These figures clearly demonstrate the the utility-optimality of the UMW+ policy.


\begin{figure}
 %\center 
 \hspace{-8pt}
 \begin{minipage}{.1\textwidth}
% \begin{overpic}[scale=0.15]{D_var_V_5_50_100_theta_1}
\begin{overpic}[scale=0.26]{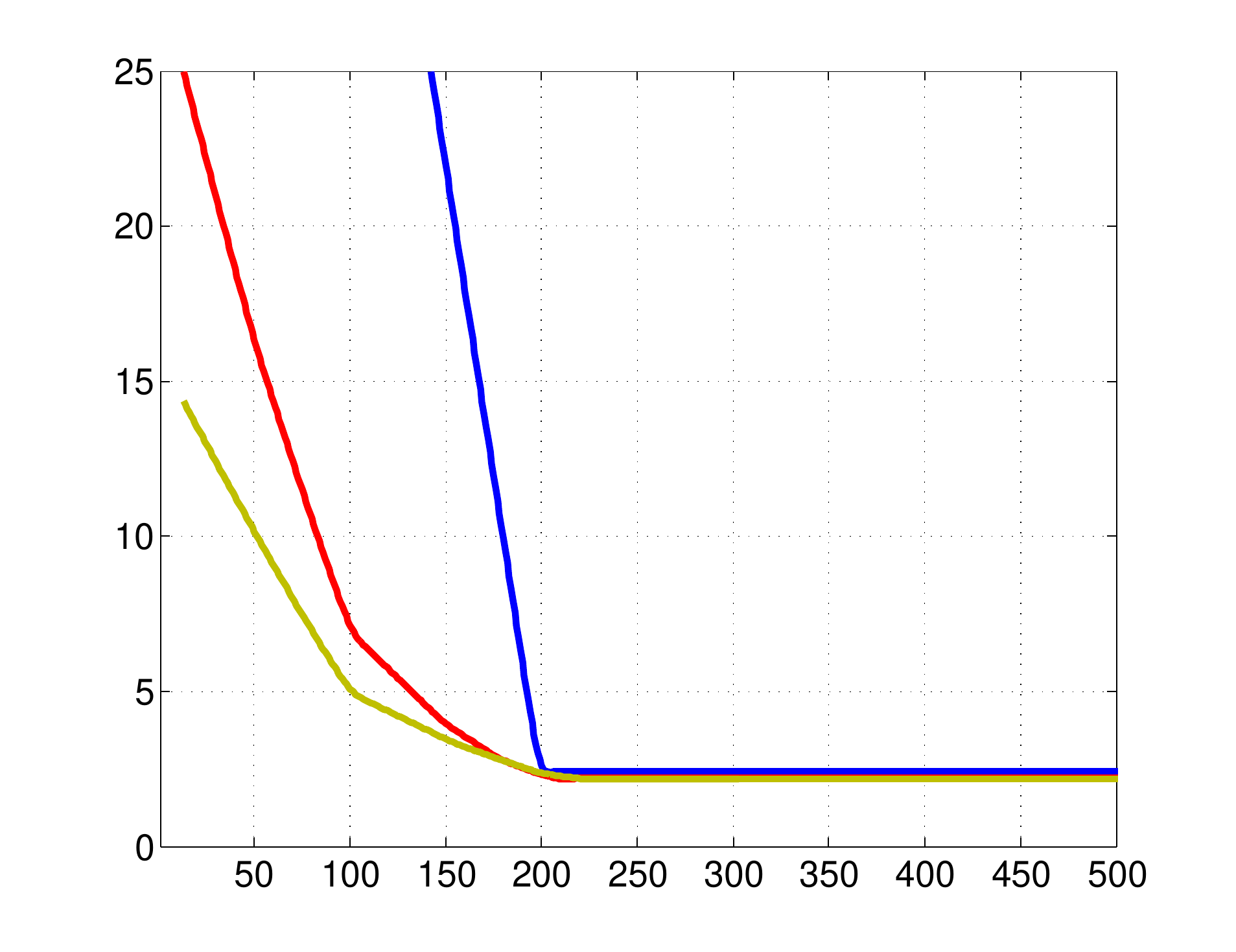}
 \put(42,32){\tiny{$V=5$}}
 \put(22,45){\tiny{$V=50$}}
 \put(15,16){\tiny{$V=100$}}
 \put(46,0.5){ \scriptsize{$t$}}
 \put(3,14){\rotatebox{90}{\scriptsize{Dual Objective $D(\bm{q}(t))$}}}
 \put(47,-3){\scriptsize{{(a)}}}
 \put(45,75){\scriptsize{$\theta=1$}}
 \end{overpic}
 \end{minipage}
 %\caption{\footnotesize{Variation of the Dual Objective with time (iterations) under Subgradient Descent for step-size $\theta=1$.}}
 \hspace{75pt}
 \begin{minipage}{.1\textwidth}
% \begin{overpic}[scale=0.23]{../NUM_Simulations/unicast_NUM/D_var_V_5_50_100_theta_pt1}
 \begin{overpic}[scale=0.26]{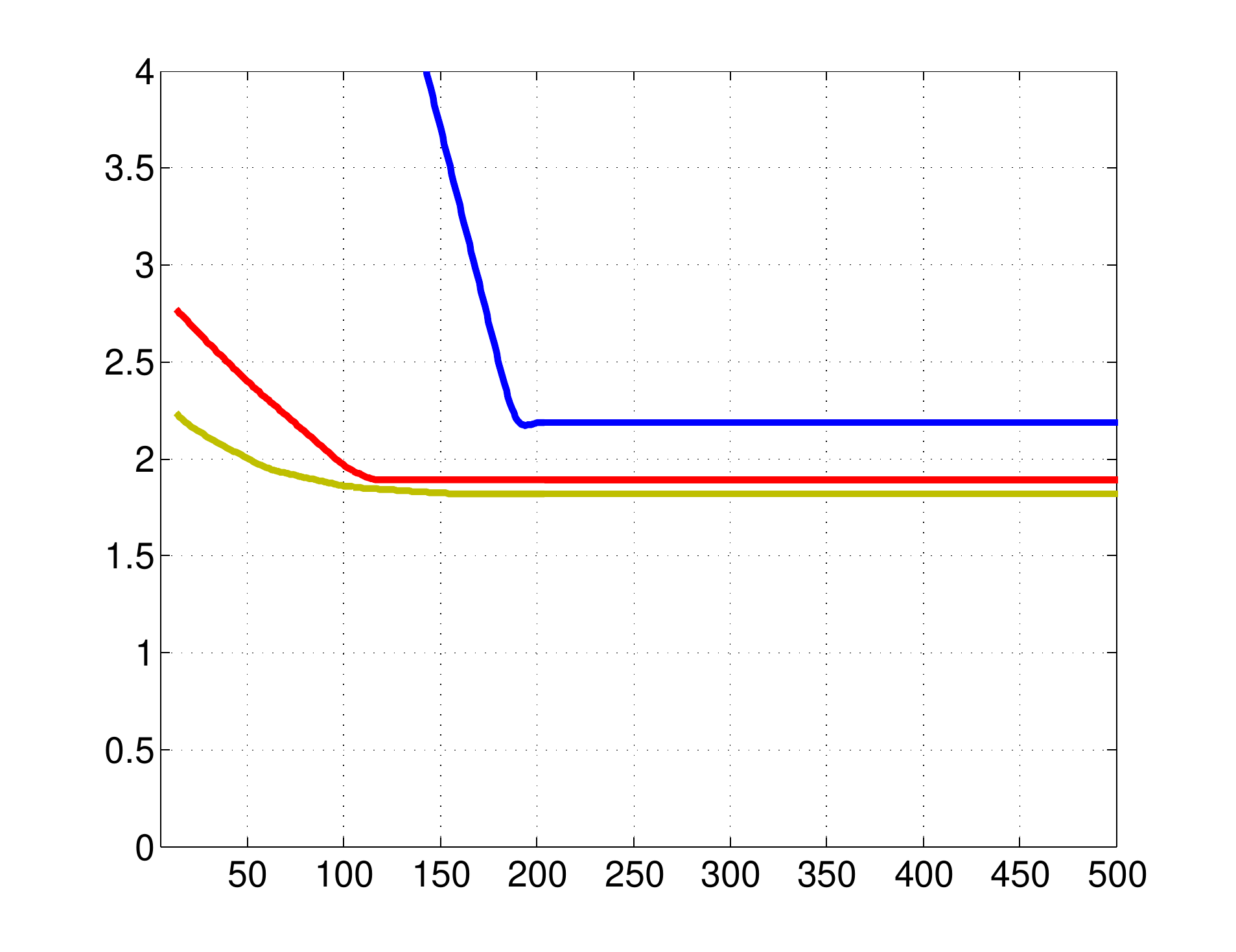}
 \put(42,48){\tiny{$V=5$}}
 \put(22,44){\tiny{$V=50$}}
 \put(16,33){\tiny{$V=100$}}
 \put(50,-1){ \scriptsize{$t$}}
 \put(51,-5){\scriptsize{{(b)}}}
 \put(1.5,14){\rotatebox{90}{\scriptsize{Dual Objective $D(\bm{q}(t))$}}}
  \put(42,75){\scriptsize{$\theta=0.1$}}
  \end{overpic}
  \end{minipage}
  \caption{\footnotesize{Variation of the Dual Objective with time (iterations) under subgradient descent for step-sizes $\theta=1$, and $\theta=0.1$.}}
  \label{time_var}
\end{figure}
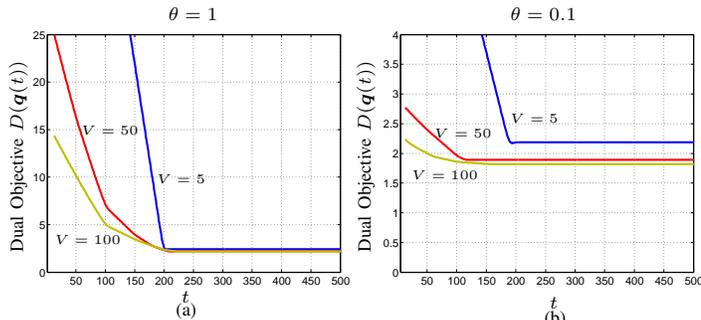
Figure \ref{time_var} shows the temporal dynamics of the subgradient algorithm for different $V$ parameters, and two different step-sizes, $\theta=1$ and $\theta=0.1$. Note that, any solution to the dual problem gives an upper bound of the primal objective (weak duality). The optimal solution of the dual minimization problem has zero duality gap with the optimal primal solution. It is evident from the plots in Figure \ref{time_var} that the speed of convergence increases and the optimality gap decreases with the increase of the parameter $V$. This observation is consistent with the statement of Theorem \ref{subgrad_convergence}, which states that the solution obtained by the dual algorithm lies within an additive gap of $O(1/V)$ from the optimal utility $U^*$.  
\begin{figure}
 \hspace{-20pt}
\begin{minipage}{.1\textwidth}
 \begin{overpic}[scale=0.13]{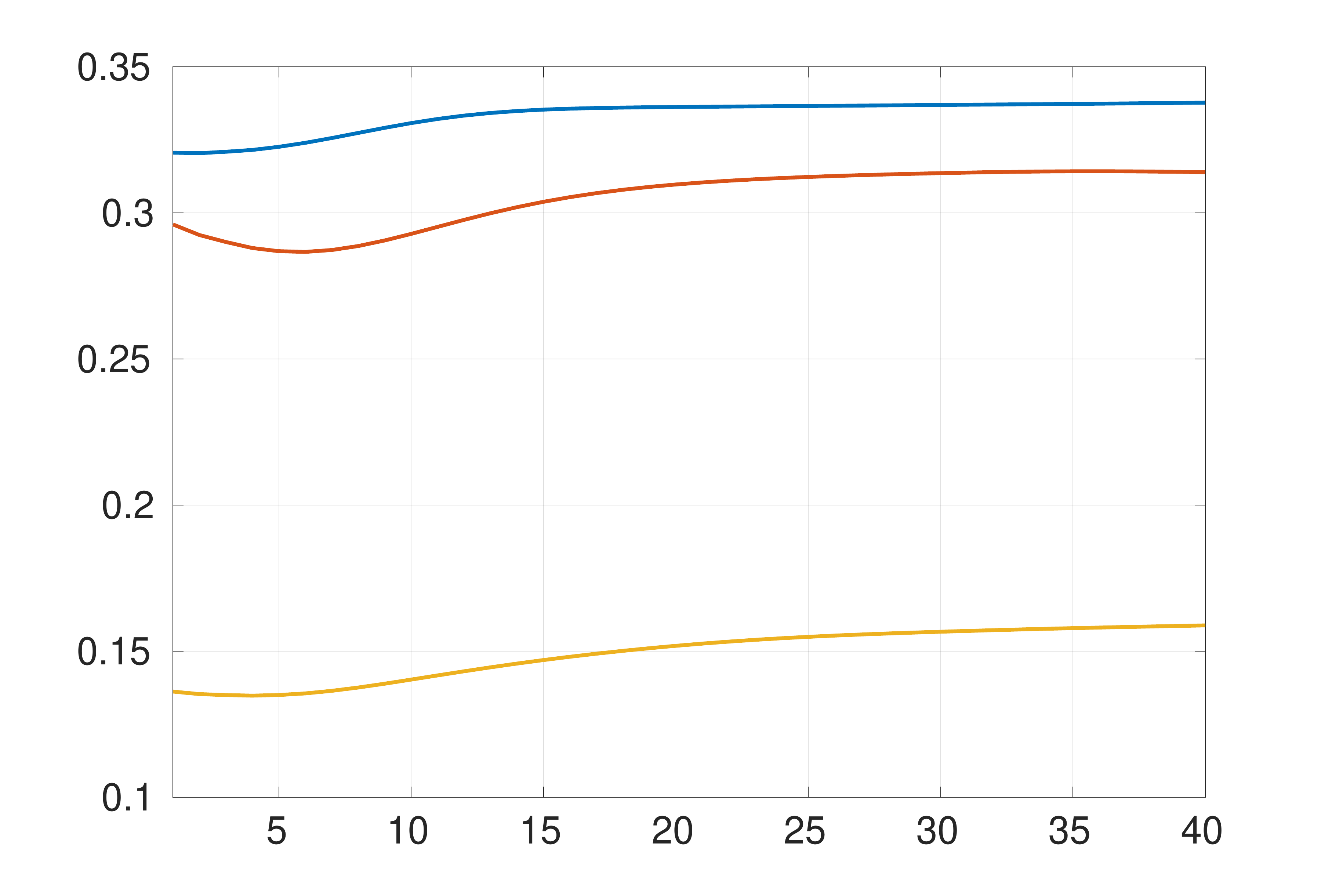}
 \put(42,56){\tiny{$\textcolor{blue}{p_{\textrm{ON}}=1.0}$}}
 \put(20,44){\tiny{$\textcolor{red}{p_{\textrm{ON}}=0.6}$}}
 \put(20,11){\tiny{$\textcolor{brown}{p_{\textrm{ON}}=0.2}$}}
 \put(50,-2){ \scriptsize{$V$}}
 \put(52,-10){\small{(a)}}
 \put(0,17){\rotatebox{90}{\scriptsize{Average Utility $U^*$}}}
  \end{overpic}
  \end{minipage}
   \hspace{90pt}
  \begin{minipage}{.1\textwidth}
 \begin{overpic}[scale=0.13]{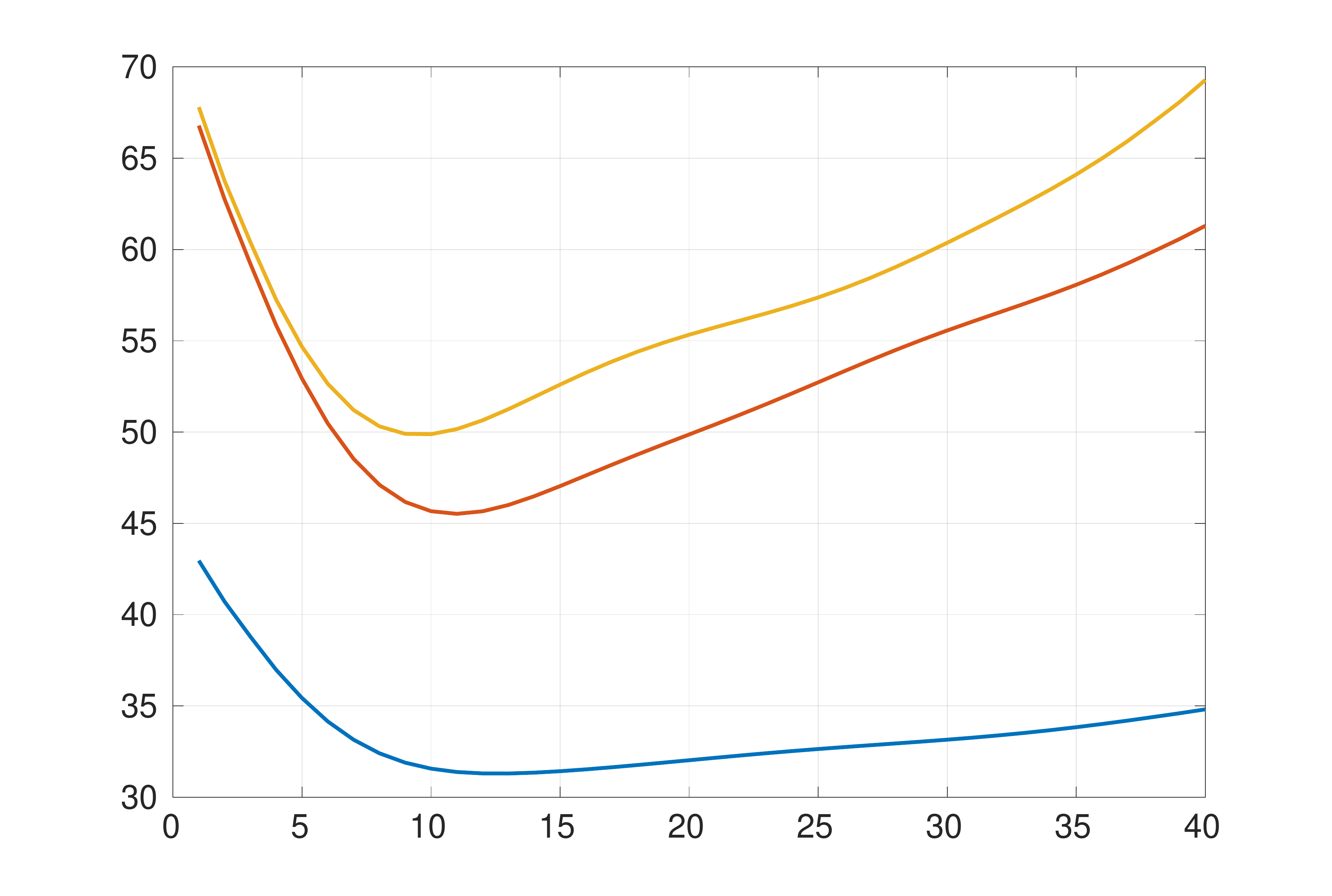}
 \put(42,12){\tiny{\textcolor{blue}{$p_{\textrm{ON}}=1.0$}}}
 \put(30,24){\tiny{\textcolor{red}{$p_{\textrm{ON}}=0.6$}}}
 \put(24,40){\tiny{\textcolor{brown}{$p_{\textrm{ON}}=0.2$}}}
 \put(50,-2){\scriptsize{ $V$}}
 \put(2,12){\rotatebox{90}{\scriptsize{Average Queue Length}}}
 \put(52,-10){\small{(b)}}
  \end{overpic}
\end{minipage}
\vspace{10pt}
\caption{\footnotesize{Variation of time-averaged Utility and total physical queue lengths with the parameter $V$ under the UMW+ policy in the broadcast setting for the time-varying grid network in Fig \ref{Network_Topolgy} (b).}}
\label{broadcast_fig}
\end{figure}
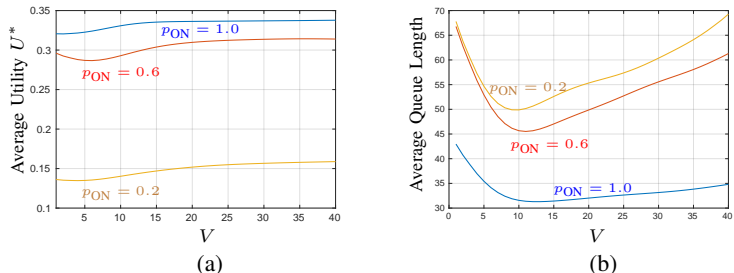
\subsection*{Broadcast Traffic in a Time-varying Wireless Network}
Next, we simulate the UMW+ policy for a time-varying wireless network, shown in Figure  \ref{Network_Topolgy} (b). Each wireless link is ON i.i.d. at every slot with probability $p_{\textrm{ON}}$. Link activations are limited by primary interference constraints. We use the same logarithmic utility function as before. Under the proposed UMW+ policy, variation of the average utility and average queue lengths with the parameter $V$ is shown in Figure \ref{broadcast_fig} (a) and \ref{broadcast_fig} (b), for three different values of $p_{\textrm{ON}}$.  As expected, with better average channel conditions (i.e., higher values of $p_{\textrm{ON}}$) higher utility is achieved with smaller in-network queue lengths.

\section{Conclusion} \label{conclusion_section}
In this paper, we have proposed the first network control policy, called UMW+, to solve the Utility Maximization Problem with multiple types of concurrent traffic, including unicast, broadcast, multicast, and anycast. The proposed policy effectively exploits the novel idea of \emph{precedence-relaxation} of a multi-hop network. We relate the UMW+ policy to the subgradient iterations of an associated dual problem. The dual objective function of the associated static NUM problem has been characterized in terms of the Fenchel conjugates of the associated utility function. Finally, illustrative simulation results have been provided for both wired and wireless networks under different input traffic settings. 
%Future work will include designing more efficient algorithms to solve the dual problem and translating these algorithms to provably optimal dynamic policies to solve the NUM problem. 

\bibliographystyle{IEEEtran}
\bibliography{MIT_broadcast_bibliography}
\section{Appendix} \label{appendix_section}

\subsection{Proof of Optimality of UMW+} \label{optimality_proof}
\begin{proof}
It is well-known that (Theorem 4.5 of \cite{neely2010stochastic}), for arbitrarily small $\epsilon>0$, there exists a stationary randomized policy ${\textrm{RAND}}$, which takes randomized actions (packet admission, routing and link scheduling) based only on the current realization of the network state $\bm{\sigma}(t)$\footnote{This class of policies have been referred to as the $\omega$-only policies in \cite{neely2010stochastic}.}, and achieves near-optimal performance: 
\begin{eqnarray}
 \mathbb{E}\sum_k U_k\big(A^{{\textrm{RAND}}}_k(t)\big) &\geq& U^* -\epsilon \label{perf1}\\
 \mathbb{E}(\tilde{A}^{{\textrm{RAND}}}_e(t) - c_e\mu^{{\textrm{RAND}}}_e(t) \sigma_e(t)) &\leq& \epsilon, \hspace{3pt}\forall e  \in E,\label{perf2}
\end{eqnarray}
where the expectations are taken over the randomness of the network configuration process $\{\bm{\sigma}(t)\}_{t\geq 0}$ and the randomized actions of the policy $\textrm{RAND}$. \\
Combining Eqns. \eqref{drift_expr} with the objective $(*)$, under the action of the \textbf{UMW+} policy, we have 
\begin{eqnarray*}
 &&\Delta(\bm{\tilde{Q}}(t), \bm{\sigma}(t))-V\sum_k U_k(A^k(t)) \leq \\
% && B+ \\
&& B+ 2\mathbb{E}\bigg(\sum_k (A^k)^{\textrm{UMW+}}(t)\big(\sum_{e} \tilde{Q}_e(t)\mathds{1}(e \in (T^k)^{\textrm{UMW+}}(t))\big)\\&&|\bm{\tilde{Q}}(t), \bm{\sigma}(t)\bigg)- 2\mathbb{E}\bigg(\sum_e \tilde{Q}_e(t)\mu^{\textrm{UMW+}}_e(t) c_e\sigma_e(t)|\bm{\tilde{Q}}(t),\\
&& \bm{\sigma}(t)\bigg)-V\sum_k U_k((A^k)^{\textrm{UMW+}}(t)) \stackrel{\textbf{(a)}}{\leq} \\
% &&\leq B+ \\
 &&B+ 2\mathbb{E}\bigg(\sum_k (A^k)^{\textrm{RAND}}(t)\big(\sum_{e} \tilde{Q}_e(t)\mathds{1}(e \in (T^k)^{\textrm{RAND}}(t))\big)\\
 &&|\bm{\tilde{Q}}(t), \bm{\sigma}(t)\bigg)- 2\mathbb{E}\bigg(\sum_e \tilde{Q}_e(t)\mu_e^{\textrm{RAND}}(t) c_e\sigma_e(t)|\bm{\tilde{Q}}(t),\\
 &&\bm{\sigma}(t)\bigg)-V\mathbb{E}\bigg(\sum_k U_k((A^k)^{\textrm{RAND}}(t))|\tilde{\bm{Q}}(t),\bm{\sigma}(t)\bigg),
\end{eqnarray*}
where the inequality (a) follows from the construction of the \textbf{UMW+} policy in Section \ref{virtual_queues}, which is defined to be minimizing the upper-bound of the conditional drift plus penalty, given in the objective $(*)$. Rearranging the RHS of the above equation, and recalling that the actions of the $\pi^\textrm{RAND}$ policy are \emph{independent} of the virtual queue lengths $\bm{\tilde{Q}}(t)$, we have 
\begin{eqnarray} \label{drift_ineq2}
 &&\Delta(\bm{\tilde{Q}}(t), \bm{\sigma}(t))-V\sum_k U_k(A^k(t)) \leq B+ \nonumber \\
 && 2\sum_{e}\bigg(\tilde{Q}_e(t)\mathbb{E}(\tilde{A}^{\textrm{RAND}}_e(t)- \mu_e^{\textrm{RAND}}(t) c_e\sigma_e(t))  \big)| \bm{\sigma}(t)\bigg) \nonumber  \\
 &&-V\mathbb{E}\bigg(\sum_k U_k((A^k)^{\textrm{RAND}}(t))|\bm{\sigma}(t)\bigg)
\end{eqnarray}
Taking expectations of the inequality \eqref{drift_ineq2} throughout w.r.t. the random network state process $\{\bm{\sigma}(t)\}_{t\geq 0}$, we have 
\begin{eqnarray*}
  &&\Delta(\bm{\tilde{Q}}(t))-V\mathbb{E}\big(\sum_k U_k(A^k(t))|\bm{\tilde{Q}}(t)\big) \leq B+ \\
   && 2\sum_{e}\bigg(\tilde{Q}_e(t)\mathbb{E}(\tilde{A}^{\textrm{RAND}}_e(t)- \mu_e^{\textrm{RAND}}(t) c_e\sigma_e(t)) \bigg) \nonumber  \\
 &&-V\mathbb{E}\bigg(\sum_k U_k((A^k)^{\textrm{RAND}}(t))\bigg)
\end{eqnarray*}
Now using the properties of the $\pi^{\textrm{RAND}}$ policy from Eqns. \eqref{perf1} and \eqref{perf2}, we have 
\begin{eqnarray*}
 &&\Delta(\bm{\tilde{Q}}(t))-V\mathbb{E}\big(\sum_k U_k(A^k(t))|\bm{\tilde{Q}}(t)\big) \leq B+ \nonumber \\
 &&2 \epsilon \sum_e \tilde{Q}_e(t)-VU^*+V\epsilon 
\end{eqnarray*}
Since the above equation holds for arbitrarily small $\epsilon >0$, taking $\epsilon \to 0$, we have the bound
\begin{eqnarray*}
 \Delta(\bm{\tilde{Q}}(t))-V\mathbb{E}\big(\sum_k U_k(A^k(t))|\bm{\tilde{Q}}(t)\big) \leq B-VU^*.
\end{eqnarray*}
Taking expectation of the above inequality w.r.t. the Virtual Queue length process $\{\bm{\tilde{Q}}(t)\}_{t \geq 0}$, we have
\begin{eqnarray} \label{key_ineq}
&& \big(\mathbb{E}L(\tilde{\bm{Q}}(t+1))- \mathbb{E}L(\tilde{\bm{Q}}(t))\big) -V\mathbb{E}\big(\sum_k U_k(A^k(t))\big) \nonumber  \\
&\leq & B-VU^*.
\end{eqnarray}
Adding up the inequalities in Eqn. \eqref{key_ineq} corresponding to the time $t=0, \ldots T-1$, telescoping the LHS, and dividing both sides by $T$, we obtain the following performance guarantee under the action of the \textbf{UMW+} policy:
\begin{eqnarray} \label{key_ineq2}
 \frac{1}{T}\mathbb{E}L(\tilde{\bm{Q}}(T))  \leq  B+\frac{1}{T}\mathbb{E}L(\tilde{\bm{Q}}(0)) + \nonumber \\
 V (\frac{1}{T}\sum_{t=1}^{T}\mathbb{E}\big(\sum_k U_k(A^k(t))\big)-U^*)
\end{eqnarray}
Since, $U^*$ is the optimal utility achievable by any feasible policy $\pi \in \Pi$, we have 
\begin{eqnarray} \label{opt_util}
 &&\limsup_{T\to \infty} \frac{1}{T}\sum_{t=1}^{T}\mathbb{E}\big(\sum_k U_k(A^k(t))\big)\nonumber \\
 &\stackrel{(a)}{\leq}& \limsup_{T\to \infty} \mathbb{E}\big(\sum_k U_k(\frac{1}{T}\sum_{t=1}^{T}A^k(t))\big) \nonumber 
 \leq  U^*,
\end{eqnarray}
where the inequality (a) follows by Jensen's inequality applied to the concave functions $U_k(\cdot), \forall k$. \\
Taking limits of both sides of the inequality \eqref{key_ineq2} and using \eqref{opt_util}, we conclude that 
\begin{eqnarray*} \label{bdd_lyap}
 \limsup_{T \to \infty}  \frac{1}{T}\mathbb{E}L(\tilde{\bm{Q}}(T))  \leq B 
\end{eqnarray*}
Since, $L(\bm{\tilde{Q}}(T))\stackrel{(\textrm{def})}{=}\sum_e \tilde{Q}_e^2(T)$, the above equation implies 
\begin{eqnarray*} 
  \limsup_{T \to \infty} \frac{1}{T}\mathbb{E}\big(\tilde{Q}_e^2(T)\big) \leq B, \hspace{10pt} \forall e \in E.
\end{eqnarray*}
By definition of $\limsup$, this implies that there exists a finite time $T^*$, such that, for all $T\geq T^*$, we have 
\begin{eqnarray*}
 \frac{1}{T}\mathbb{E}\big(\tilde{Q}_e^2(T)\big) \leq B+1, 
\end{eqnarray*}
i.e., $ \mathbb{E}\big(\tilde{Q}_e^2(T)\big) \leq T(B+1)$. Applying Jensen's inequality \cite{durrett2010probability}, we conclude 
\begin{eqnarray}\label{bd2}
 \mathbb{E}\big(\tilde{Q}_e(T)\big) \leq \sqrt{T}\sqrt{B+1}, \hspace{10pt} \forall T\geq T^*.
\end{eqnarray}

Using the bounded admission assumption in Section \ref{admissible}, the total arrival to each virtual queue $\tilde{Q}_e$ is upper bounded by $|\mathcal{K}|A_{\max}$. Hence, we have the following trivial deterministic upper bound on the virtual queue sizes at time $T$:
\begin{eqnarray} \label{bd1}
 \tilde{Q}_e(T) \leq |\mathcal{K}|A_{\max}T, \hspace{10pt} \forall e  \in E, \textrm{ a.s.}
 \end{eqnarray}
 Thus, the sequence of non-negative random variables $\{\frac{\tilde{Q}_e(T)}{T}\}_{T\geq 1}$ are \emph{uniformly bounded}. 
 
 Using the Bounded Convergence Theorem \cite{durrett2010probability} \footnote{ We are implicitly \emph{assuming} that the limiting random variable is well-defined almost surely.}, we conclude that
\newpage
 \begin{eqnarray} \label{as_conv}
  \mathbb{E}\bigg(\lim_{T \to \infty} \frac{\tilde{Q}_e(T)}{T} \bigg) = \lim_{T \to \infty} \frac{\mathbb{E}(\tilde{Q}_e(T))}{T}\stackrel{(b)}{=}0,
 \end{eqnarray} 
where the equality (b) follows from the bound in Eqn. \eqref{bd2}. Since the random variables $\{\tilde{Q}_e(T)\}_{T\geq 1}$ are non-negative, Eqn. \eqref{as_conv} implies that
\begin{eqnarray*}
 \lim_{T\to \infty} \frac{\tilde{Q}_e(T)}{T} =0, \hspace{10pt} \textrm{a.s.}, \forall e \in E.
\end{eqnarray*}
This proves part (a) of Theorem \eqref{VQ_stability}. To establish part (b), \emph{i.e.},  utility optimality of the \textbf{UMW+} policy, we consider Eqn. \eqref{key_ineq2} again and using non-negativity of the Lyapunov function $L(\bm{\tilde{Q}}(T))$, we have 
\begin{eqnarray} \label{U-o}
&&\mathbb{E}\big(\sum_{k}U_k(\frac{1}{T}\sum_{t=1}^{T}A_k(t)))\nonumber \\
&\stackrel{(a)}{\geq}& 
 \frac{1}{T}\sum_{t=1}^{T}\mathbb{E}\big(\sum_k U_k(A^k(t))\big) \geq U^* - \frac{B}{V},
\end{eqnarray}
where the inequality (a) follows from Jensen's inequality applied to the concave functions $U_k(\cdot), \forall k$.
Utility optimality of \textbf{UMW+} follows by taking the limit $T \to \infty$ of Eqn. \eqref{U-o}. 
\end{proof}

\end{document}